\documentstyle{article}

\newtheorem{theorem}{Theorem}

\newtheorem{definition}[theorem]{Definition}
\newtheorem{example}[theorem]{Example}
\newtheorem{remark}[theorem]{Remark}

\def\QED{\quad\blackslug\lower 8.5pt\null}

\begin{document}

\begin{center}
{\Large \bf  LIGHTLIKE HYPERSURFACES} 

\vspace*{2mm}

{\Large \bf ON MANIFOLDS ENDOWED}
 
\vspace*{2mm}
{\Large \bf  WITH A CONFORMAL STRUCTURE}
 
\vspace*{2mm}

{\Large \bf  OF LORENTZIAN SIGNATURE}\footnote{{\bf 
1991 MS classification}: 53A30, 53B25.

\hspace*{1mm} {\bf Keywords and phrases:} 
Pseudoconformal structure,  Lorentzian signature, lightlike 
hypersurface,   isotropic geodesics, singular point, 
invariant normalization, affine connection.}

\vspace*{3mm}
{\large M.A. Akivis and  V.V. Goldberg}

\end{center}

\vspace*{5mm}

{\footnotesize{\it Abstract}. The authors study the geometry of 
lightlike hypersurfaces on manifolds $(M, c)$ endowed with 
a pseudoconformal structure $c = CO (n - 1, 1)$ of Lorentzian 
signature. Such hypersurfaces are of interest in general 
relativity since they can be models of different types of 
physical horizons. On a lightlike hypersurface, the 
authors consider the fibration of isotropic geodesics 
and investigate their singular points and singular submanifolds. 
They construct a conformally invariant normalization of 
a lightlike hypersurface intrinsically connected with 
its geometry and investigate affine connections induced by 
this normalization. The authors also consider  special 
classes of lightlike hypersurfaces. In particular, they 
investigate lightlike hypersurfaces for which the elements 
of the constructed normalization are integrable.}
\vspace*{5mm}

\setcounter{equation}{0}

\setcounter{section}{-1}

\section{Introduction}  

The pseudo-Riemannian manifolds $(M, g)$ of Lorentzian signature 
play a special role in geometry and physics: they generate 
models of spacetime of general relativity. In the tangent 
space $T_x$ at a point $x$ of such a manifold, a real isotropic 
cone $C_x$ is  invariantly defined. From physical 
point of view, this cone is the light cone---the trajectories of 
light impulses emanating from the point $x$ are tangent to the cone $C_x$.

A hypersurface $V^{n-1}$ on an $n$-dimensional manifold $M$ 
 of  Lorentzian signature that is tangent to the cone 
$C_x$ at each  point $x \in V$ is called lightlike. The lightlike 
hypersurfaces are also of interest for general relativity since 
they produce models of different type of horizons 
(event horizons, Cauchy's horizons,  Kruskal's horizons---
see, for example,  [Ch 83] and [MTW 73]). Lightlike hypersurfaces 
are also studied in the theory of electromagnetism. This is the 
reason that there are many papers and two recent books 
[DB 96] and [Ku 96] 
in which lightlike hypersurfaces are investigated.

Many events and objects of general relativity are invariant under 
conformal transformations of a metric (see [AG 96], Ch. 4 and 
Ch. 5). 
In particular, a lightlike hypersurface is an example of 
the  objects that are invariant under 
conformal transformations of a metric. Hence it is appropriate to study lightlike hypersurfaces 
not only on a pseudo-Riemannian manifold $(M, g)$ of Lorentzian 
signature but also on a manifold endowed with a conformal 
structure of Lorentzian signature. 

In the present paper we 
study lightlike hypersurfaces on a differentiable manifold $M^n$ 
endowed with a pseudoconformal structure $CO (n-1, 1)$ of 
Lorentzian signature. We will denote such manifolds by $(M, c)$ 
where $\dim \, M = n$ and $c = CO (n-1, 1)$ is a conformal 
structure of signature $(n-1, 1)$.

Let us describe the contents of the paper. In Section {\bf 1} 
we write the basic equations of the manifold $(M, c)$ and 
consider pseudoconformal spaces $(C^n_1)_x$ of Lorentzian 
signature tangent to the manifold $(M, c)$ at its point $x$. 
It appeared that it is very convenient to use the 
Darboux representation of a space $(C^n_1)_x$ as a hyperquadric 
$(Q^n_1)_x$ of a projective space $P^{n+1}$.

As we proved in [AG 96] (see also [AG 97]), 
the isotropic geodesics of 
a pseudo-Riemannian manifold $(M, g)$ are invariant under 
conformal transformations of a metric. Thus they can be 
considered on a manifold $(M, c)$ endowed with a pseudoconformal 
structure.  Under the development of 
 the manifold $(M, c)$  onto a 
hyperquadric $(Q_1^n)_x$ along an  isotropic 
geodesic, the latter is mapped into a rectilinear generator of 
the hyperquadric $(Q_1^n)_x$.

In Section {\bf 2} we consider differential geometry of 
lightlike hypersurfaces $V^{n-1} \subset (M, c)$. In this section 
we construct a first-order frame bundle associated with 
$V^{n-1}$; define a screen distribution $S$ (see [DB 96]) and 
a field $N$ of normalizing isotropic straight lines 
that is conjugate to $S$; write the basic equations of 
lightlike hypersurfaces and prove the existence theorem 
for such hypersurfaces; prove that a lightlike hypersurface 
carries $(n-1)$-parameter family of isotropic geodesics 
each of which possesses  $n - 2$ real singular points 
if each of them is counted as many times 
as its multiplicity; and prove that 
under the development of  the hypersurface  $V^{n-1}$  onto a 
hyperquadric $(Q_1^n)_x$ along an  isotropic 
geodesic $l$, to $V^{n-1}$ there corresponds a ruled 
hypersurface which has the same tangent subspace at 
all regular points of its rectilinear generator. 

In Section {\bf 3} we introduce the basic geometric objects 
and tensors defined in a second-order neighborhood 
of a point of a lightlike hypersurface $V^{n-1}$ as well as 
geometric images associated with these objects and tensors. 
In particular, on each isotropic geodesic $l$ of the hypersurface 
$V^{n-1}$ we construct the harmonic pole $H$ of 
the point $x \in l \subset V^{n-1}$ with respect to 
singular points of the generator \nolinebreak $l$.

In studying submanifolds on a manifold $M$ endowed with 
a differential geometric structure defined by a group $G$, 
one of the most important problems is a construction of 
an invariant normalization and an affine connection intrinsically 
connected with the geometry of a submanifold in question 
(see [La 53], [No 76], [AG 93, 96]). In some simple cases 
such a normalization and a connection are defined 
in a first-order neighborhood. This is the case for 
submanifolds of the Riemannian manifold and for 
spacelike and timelike submanifolds of the pseudo-Riemannian 
manifold. However, for lightlike hypersurfaces 
of a  pseudo-Riemannian manifold $(M, g)$ as well as 
of a manifold $(M, c)$ endowed with a pseudoconformal 
structure, for constructing such a normalization 
and an affine connection elements of higher order 
differential neighborhoods are needed. 

In Section {\bf 4} we give a geometric construction (called 
normalization) defining an affine connection 
on $V^{n-1} \subset (M, c)$. However, the main purpose of this 
paper is a construction of an invariant normalization of 
$V^{n-1}$ and an affine connection induced by this 
normalization. A solution of this problem is presented 
in Sections ${\bf 5}$ and ${\bf 6}$. In Section ${\bf 5}$ 
we construct geometric objects defined 
in a third-order neighborhood of $x \in V^{n-1}$ and apply them 
to construct a screen distribution $S$, whose elements are 
subspaces $S_x \subset T_x (V^{n-1})$ passing through the point 
$x$, and a complementary screen distribution $\widetilde{S}$,  
whose elements are 
subspaces $\widetilde{S}_H$ passing through the harmonic pole 
$H$ of the point $x$. The above construction can be done 
provided that two affinors are nondegenerate on $V^{n-1}$. 
One of these affinors is defined in a second-order neighborhood 
of $x \in V^{n-1}$, and the second one in a third-order 
neighborhood of $x \in V^{n-1}$.  

Under the same assumptions in Section {\bf 6} we construct 
a one-component geometric object which is defined in a 
fourth-order neighborhood of $x \in V^{n-1}$ and determines 
a point $C_n$ on a normalizing isotropic straight line. 
All these geometric objects are intrinsically defined by the 
geometry of $V^{n-1}$. A geometric meaning of the geometric 
images associated with the constructed objects is also found. 
The invariant normalization of $V^{n-1}$ we have constructed 
induces a torsion-free affine connection $\gamma_1$ whose 
fundamental group is the group $G_1 = GL (n-1, {\bf R})$. 
The curvature tensor of this connection is expressed in terms of 
elements of a fifth-order neighborhood of $x \in V^{n-1}$.

In Section {\bf 7} we investigate the problem of integrability of 
the screen distributions $S$ and $\widetilde{S}$. We prove that 
they 
are simultaneously integrable or simultaneously nonintegrable and 
find conditions of their integrability. If the 
 distributions $S$ and $\widetilde{S}$ are integrable, then the 
 congruence of normalizing isotropic straight lines 
$\widetilde{l}$ is stratified into a one-parameter family of 
lightlike hypersurfaces.

Finally in Section {\bf 7}, in addition to 
the torsion-free affine connection $\gamma_1$ induced by 
the constructed invariant normalization, we find 
another affine connection $\gamma_2$ whose fundamental group 
$G_2$ is isomorphic to the group 
${\bf R}^+ \times GL (n-2, {\bf R})$. The connection $\gamma_2$ 
is not torsion-free and  defined by elements a 
third-order neighborhood of $x \in V^{n-1}$. The torsion 
and curvature 
tensors of this connection are expressed in terms of 
elements of a fourth-order neighborhood of $x \in V^{n-1}$.

Note that the problem of construction of an invariant 
normalization and an invariant affine connection for 
a lightlike hypersurface $V^{n-1}$ of a pseudo-Riemannian 
manifold $(M, g)$ of Lorentzian signature (for definition 
see [ON 83]) was 
considered by many authors (see [DB 96], Ch. 4). 
However, as far as we know, an invariant 
normalization and an affine connection intrinsically 
connected with the geometry of $V^{n-1}$  
were not considered in these papers. 

As to a construction of invariant normalization and an 
invariant affine connection for a lightlike hypersurface 
$V^{n-1}$ on a manifold $(M, c)$ endowed with a pseudoconformal 
structure of Lorentzian signature, such a construction  is given 
in the present paper for the first time. 

The methods developed in the present paper and the results 
obtained in it can be used to study lightlike 
hypersurfaces in the pseudo-Riemannian space $(M, g)$ of 
Lorentzian signature, the pseudoconformal space $C^n_1$, and the 
Minkowski space $R^n_1$.

The present paper is related to our papers 
[AG 98a, b, c].

\section{Basic equations of a manifold endowed with a conformal 
structure of Lorentzian signature}

{\bf 1.} A pseudoconformal structure 
$CO (n-1, 1)$ on a manifold $M$ of dimension $n$ is a set of 
conformally equivalent pseudo-Riemannian metrics with 
the same signature $(n-1, 1)$. Such a structure is called 
{\em conformally Lorentzian}.

A  metric $g$ can be given on 
$M$ by means of a nondegenerate quadratic form 
$$
g = g_{ij} du^i du^j,
$$
where $u^i, \;
i = 1, \ldots , n$, are curvilinear coordinates on  
$M$, and $g_{ij}$ are the components of the metric tensor $g$.
 On a  conformal structure the quadratic form $g$ is 
relatively invariant.

The equation $g = 0$ defines in the tangent space 
$T_x (M)$ 
a cone $C_x$ of second order called the {\em isotropic 
cone}. Thus the conformal structure can be given on 
the  manifold $M$ by a field of cones of  second order.

The cone $C_x \subset T_x (M)$ remains invariant under 
transformations of the group 
$$
 G = {\bf SO} (n-1, 1)  \times {\bf H}, 
$$
where ${\bf SO} (n-1, 1)$ is the special $n$-dimensional 
pseudoorthogonal 
group of signature $(n-1, 1)$ (the connected component of the 
unity of the pseudoorthogonal group ${\bf O} (n-1, 1)$), 
and ${\bf H}$ is the group of homotheties. Note that the group 
$G$ acts in the tangent space $T_x (M)$. 
Thus  the  pseudoconformal structure $CO (n-1, 1)$ is a 
$G$-structure defined on the manifold $M$ by the group $G$ 
indicated above. For  the pseudoconformal structure $CO (n-1, 1)$ 
the isotropic cone is real. Note that pseudo-Riemannian 
structures of arbitrary signature were studied in the book 
[ON 83] (they were called there "semi-Riemannian").

{\bf 2.} We consider the manifold $M$, associate with any point 
$x \in M$ its tangent space $T_x (M)$, and define the frame 
bundle whose base is the manifold $M$ and the fibers are the 
families of vectorial frames $\{e_1, \ldots , e_n\}$ in $T_x (M)$ 
defined up to a transformation of the general linear group 
${\bf GL} (n)$. The frames indicated above are called the {\em 
frames of first order}. They form the first-order 
frame bundle which 
we will denote by ${\cal R}^1 (M)$. Let us denote by $\{\omega^1, 
\ldots , \omega^n\}$ the co-frame dual to the frame 
$\{e_1, \ldots , e_n\}$: 
$$
\omega^i (e_j) = \delta^i_j.
$$ 
Then an arbitrary vector $\xi \in T_x (M)$ can be written as 
$$
\xi = \omega^i (\xi) e_i.
$$
The forms $\omega^i$ can be considered as differential forms on 
the manifold $M$ if we assume that $\xi = dx$ is the differential 
of the point $x \in M$. Thus the quadratic form $g$ can be 
written as 
\begin{equation}\label{eq:1}
g = g_{ij} \omega^i \omega^j.
\end{equation}

{\bf 3.} 
The structure equations of the $CO (n-1, 1)$-structure 
can be reduced to the following form (see [AG 96], Section 
{\bf 4.1}):
\begin{equation}\label{eq:2}
d \omega^i = \omega_0^0 \wedge \omega^i 
+ \omega^j \wedge \omega^i_j, 
\end{equation}
\begin{equation}\label{eq:3}
 d \omega_0^0 = \omega^i \wedge \omega^0_i, 
\end{equation}
\begin{equation}\label{eq:4}
d \omega^i_j = \omega_j^0 \wedge \omega^i + 
 \omega^k_j \wedge \omega_k^i + g_{jk} \omega^k \wedge 
g^{il} \omega^0_l + C^i_{jkl} \omega^k \wedge \omega^l,
\end{equation} 
\begin{equation}\label{eq:5}
d \omega_i^0 = \omega^0_i \wedge \omega_0^0
 + \omega_i^j \wedge \omega^0_j 
+ C_{ijk}  \omega^j  \wedge  \omega^k. 
\end{equation}
and the metric tensor $g_{ij}$ satisfies the equations
\begin{equation}\label{eq:6}
d g_{ij} - g_{ik} \omega_j^k  -  g_{kj} \omega_i^k  = 0.
\end{equation}
Note that in equations (2)--(6) the forms $\omega^i$ 
are defined in a first-order frame bundle, 
the 1-forms $\omega^i_j$ and a scalar 
 1-form $\omega_0^0$ in a second-order frame bundle, and 
a covector form $\omega^0_i$  in the third-order 
frame bundle. 

For $C^i_{jkl} =  C_{ijk} = 0$,  equations 
(2)--(6) coincide with the structure equations 
 the pseudoconformal space $C^n_1$. 
For this  reason  the object $\{C^i_{jkl},  C_{ijk}\}$ is 
called the 
{\em curvature object of the pseudoconformal 
structure $CO (n-1, 1)$}.

 The quantities $C^i_{jkl}$ form a $(1, 3)$-tensor which 
is called the {\em Weyl tensor} or the {\em tensor of conformal 
curvature} of the structure $CO (n-1, 1)$. 

Consider also the covariant tensor of conformal curvature
\begin{equation}\label{eq:7}
C_{ijkl} = g_{im}  C^m_{jkl}.
\end{equation}
This tensor allows us to write relations between the components 
of the tensor of conformal curvature in more 
convenient form:
\begin{equation}\label{eq:8}
\renewcommand{\arraystretch}{1.3}
\left\{
\begin{array}{ll}
C_{ijkl} = - C_{jikl} = - C_{ijlk},  & C_{ijkl} = C_{klij}, \\ 
C_{ijkl} + C_{iklj} + C_{iljk} = 0.
\end{array} 
\right.
\renewcommand{\arraystretch}{1}
\end{equation}
In addition, the tensor $C^i_{jkl}$ is trace-free:
\begin{equation}\label{eq:9}
C^i_{jki} = 0.
\end{equation}
(see, for example, [AG 96], Section {\bf 4.1}). 

The quantities $C_{ijk}$, that do not form 
a tensor, satisfy the conditions 
\begin{equation}\label{eq:10}
C_{ijk} = - C_{ikj}.
\end{equation}

Note also that the tensor of conformal curvature 
$C^i_{jkl}$ is defined in a third-order neighborhood 
of the structure $CO (n - 1, 1)$, and the quantities 
$C_{ijk}$ are defined in its fourth-order neighborhood,
 The $CO (n-1, 1)$-structure itself is a $G$-structure of 
finite type two 
 (see [AG 96], Section {\bf 4.1}). 
For $n \geq 4$, the condition $C^i_{jkl} = 0$  is 
necessary and sufficient for a manifold $(M, c)$ to be 
conformally flat (see [AG 96], Section {\bf 4.1}).

{\bf 4}. The 1-forms $\omega^i = \omega^i _0, \;\omega^0_0, \;
\omega^i_j$, and $\omega^0_i$ defined by the  
$CO (n-1, 1)$-structure on the manifold $M$ can be taken as 
components of infinitesimal displacement of a frame in 
the pseudoconformal space $C_1^n$. A conformal frame consists 
of two points $A_0$ and $A_{n+1}$ and $n$ hyperspheres $A_i$ 
passing through $A_0$ and $A_{n+1}$. The scalar products of 
the elements of this frame can be written as

\begin{equation}\label{eq:11}
\renewcommand{\arraystretch}{1.3}
\left\{
\begin{array}{ll}
(A_0, A_0) = (A_{n+1}, A_{n+1}) = 0, \;\; (A_0, A_{n+1}) = -1, \\
(A_0, A_i) = (A_{n+1}, A_i) = 0, \;\; (A_i, A_j) = g_{ij}
\end{array} 
\right.
\renewcommand{\arraystretch}{1}
\end{equation}
The equations of infinitesimal displacement 
of this frame have the form
\begin{equation}\label{eq:12}
\renewcommand{\arraystretch}{1.3}
\left\{
\begin{array}{ll}
dA_0 = \omega_0^0  A_0 & \!\!\!\! +  \omega_0^i  A_i, \\
dA_i = \omega_i^0  A_0 & \!\!\!\!  +  \omega_i^j  A_j + \omega_i^{n+1} A_{n+1}, \\
dA_{n+1} = &  \!\!\!\!  \omega_{n+1}^i  A_i +  \omega^{n+1}_{n+1} A_{n+1}, 
\end{array} 
\right.
\renewcommand{\arraystretch}{1}
\end{equation}
where 
$$
 \omega_i^{n+1} = g_{ij}  \omega_0^j, \;\; 
 \omega^i_{n+1} = g^{ij}  \omega_j^0, \;\; 
 \omega_{n+1}^{n+1} = -  \omega_0^0,
$$
and $g^{ij}$ is the inverse tensor of the tensor 
$g_{ij}$. In addition, the forms $ \omega_j^i$ 
satisfy 
the system of 
equations (6). The family of frames 
in question forms a  bundle of first-order conformal 
frames associated with 
the pseudoconformal structure $(M, c)$.

Equations (12) are completely integrable if and only if 
the  tensor of conformal curvature of 
the $CO (n-1, 1)$-structure vanishes. Then these equations 
define a fiber bundle in the whole space $C^n_1$. 
If the  tensor of conformal curvature does not vanish, the system 
(12) can be integrated along any smooth curve $x = x(t)$ 
belonging to the manifold $M$. A solution of this system 
defines a {\em development} of this line and the frame bundle 
along it on  the conformal space $C^n_1$. Moreover, if 
$x_1$ and $x_2$ are two points of the manifold $M$, 
and $l_1$ and $l_2$ are two smooth curves joining 
these points, then under integration of equations (12) along 
these curves, for the same initial conditions at the point 
$x_1$,  we obtain different results at the point $x_2$. 
The difference of these two results is defined 
by the curvature of the pseudoconformal structure 
$CO (n-1, 1)$ (see [Car 23] and also [AG 96]).

For study of conformal structures it is convenient 
to use {\em Darboux mapping}. Under the Darboux mapping 
to the conformal space $C^n_1$ there corresponds 
a hyperquadric $Q^n_1$ of  a projective space $P^{n+1}$;
to the points $A_0$ and $A_{n+1}$ there correspond 
points of the  hyperquadric $Q^n_1$ not belonging to 
a rectilinear generator of $Q_1^n$; and 
to the hyperspheres $A_i$  there correspond 
points of the space $P^{n+1}$ belonging to 
the  intersection of the  hyperplanes $T_x (Q^n_1)$ 
and $T_y (Q_1^n)$ tangent to the hyperquadric $Q_1^n$ 
at its points $x = A_0$ and $y = A_{n+1}$ (see Figure 1). 
We will 
denote the elements of a projective frame by 
the same letters which 
we used for the corresponding elements 
of a conformal frame. The equations of infinitesimal 
displacement of the projective frame in question have 
the same form (12) as the equations of infinitesimal 
displacement of the corresponding conformal frame.

\vspace*{50mm}

\begin{center}
Figure 1
\end{center}

The equation of the hyperquadric $Q^n_1$ with respect to 
the  projective frames in question has the form:
\begin{equation}\label{eq:13}
(x, x) = g_{ij} x^i x^j - 2 x^0 x^n = 0.
\end{equation}
The quadratic form $g = g_{ij} x^i x^j$ is of signature 
$(n - 1, 1)$, and the equation $g_{ij} x^i x^j = 0$ 
defines the isotropic cone $C_x$ with the vertex at $x = A_0$ on 
the hyperquadric $Q_1^n$. This cone carries 
an $(n-2)$-parameter family of rectilinear generators 
corresponding to the isotropic lines of the space $C^n_1$. 

For $\omega^i = 0$, equations (12) determine an admissible 
transformation of frames in the pseudoconformal space $(C_1^n)_x$ 
that is tangent to a manifold $(M, c)$ 
endowed a pseudoconformal structure $CO (n - 1, 1)$ 
at a point $x$.

\section{Geometry of lightlike hypersurfaces of a manifold 
endowed with a conformal structure of Lorentzian signature}

{\bf 1.} In this paper we consider 
a lightlike hypersurface $V^{n-1}$ on a manifold 
$M$ of dimension $n \geq 4$ endowed with 
a $CO (n-1, 1)$-structure 
of Lorentzian signature $(n-1, 1)$. 
A {\em lightlike hypersurface} $V^{n-1}$ 
on such a manifold is a hypersurface which is tangent to the 
isotropic cone $C_x$ at each point $x \in V^{n-1}$. 

Let $T_x (V^{n-1})$ be a tangent subspace to 
$V^{n-1}$ at a point $x$. 
In $T_x (M)$ we choose a projective frame 
such that $x = A_0$; the point $A_1$ 
belongs to the isotropic  generator of the  cone $C_x$ along 
which the subspace $T_x (V^{n-1})$ is tangent to $C_x$; 
the point $A_n$ also belongs to a rectilinear 
generator of the cone $C_x$ that does not belong to the subspace 
$T_x (V^{n-1})$; and we place the points $A_a, a = 2, \ldots , 
n - 1$,   into the $(n-2)$-dimensional  intersection of the 
subspace $T_x (V^{n-1})$ and the subspace $T_{A_0 A_n} (C_x)$ 
tangent to $C_x$ along $A_0 A_n$. Then the scalar products 
of these points can be written as
\begin{equation}\label{eq:14}
\renewcommand{\arraystretch}{1.3}
\left\{
\begin{array}{ll}
(A_1, A_1) = (A_n, A_n) = 0, \;\; (A_1, A_a) = (A_n, A_a) = 0, \\ 
(A_a, A_b) = g_{ab}, \;\; (A_1, A_n) = -1,
\end{array} 
\right.
\renewcommand{\arraystretch}{1}
\end{equation}
where $a, b = 2, \ldots , n -1$. 
The last relation in (14) is a 
result of an appropriate normalization of the 
points $A_0$ and $A_n$. The frames we have 
constructed are first-order frames associated with a 
lightlike hypersurface $V^{n-1}$ (see Figure 2).

\vspace*{65mm}

\begin{center}
Figure 2
\end{center}

The isotropic straight lines $N_x = A_0 \wedge A_n$ are called 
the {\em normalizing straight lines} of a lightlike 
hypersurface $V^{n-1}$, and the subspaces 
$S_x = A_0 \wedge A_2 \wedge \ldots \wedge A_{n-1}$ 
belonging to $T_x (V^{n-1})$ are called the 
{\em screen subspaces} of $V^{n-1}$. {\em There exists 
a bijective correspondence between the fields $N$ of 
normalizing straight lines $N_x$ of a lightlike hypersurface 
$V^{n-1}$ and its screen distributions $S$ of screen subspaces 
$S_x$.}

A normalizing field $N$ on a  lightlike hypersurface 
$V^{n-1}$ can be given with a big arbitrariness. One of 
the goals of the present paper is to find a method 
of construction of a normal field $N$, and along with 
this field also a screen distribution $S$ both intrinsically 
connected with the geometry of a  lightlike hypersurface 
$V^{n-1}$.

With respect to the  projective moving frame  
chosen in the tangent space $T_x (M)$, 
 the fundamental form $g$ of $M$ 
has the expression
\begin{equation}\label{eq:15}
g = g_{ab} \omega^a \omega^b - 2 \omega^1 \omega^n, 
\;\;\;\;\; a, b = 2, \ldots , n -1, 
\end{equation}
and the isotropic cone $C_x$ is 
determined by the equation $g = 0$.

Now the components $g_{ij}$ of the  tensor  $g$ 
are the entries of the following matrix:
\begin{equation}\label{eq:16}
(g_{ij}) = \pmatrix{0 & 0 & -1 \cr 
                                0 & g_{ab} & 0\cr 
                                -1 & 0 & 0},
\end{equation}
where $(g_{ab})$ is a nondegenerate positive definite 
matrix. 
Equations (6) and (16) imply that 
\begin{equation}\label{eq:17} 
\renewcommand{\arraystretch}{1.3}
\left\{
\begin{array}{ll}
\omega_1^n = \omega_n^1 = 0, & \omega_1^1 = - \omega_n^n, \\
 \omega_a^1 =  g_{ab} \omega^b_n, &  \omega_a^n =  g_{ab} \omega^b_1,
\\ 
d g_{ab} - g_{ac} \omega_b^c  -  g_{cb} \omega_a^c  = 0. &
\end{array} 
\right.
\renewcommand{\arraystretch}{1}
\end{equation}

Since  the points $A_1$ and $A_a$ 
of the frame  $\{A_0, A_1, A_a, A_n\}$ of $T_x (V^{n-1})$ 
belong to 
the tangent subspace $T_x (V^{n-1})$,   we have 
$$
dA_0 = \omega_0^0 A_0 + \omega^1_0 A_1 + \omega_0^a A_a. 
$$
This means that the hypersurface $V^{n-1}$ is defined 
by the following Pfaffian equation
\begin{equation}\label{eq:18}
\omega_0^n = 0, 
\end{equation}
and the forms $\omega^1$ and $\omega^a,\; a = 2, \ldots , n - 1,$
 are basis forms of the hypersurface $V^{n-1}$. 

The quadratic form $\widetilde{g}$ defining 
the conformal structure on $V^{n-1}$ has the form
$$
\widetilde{g} = g_{ab} \omega^a \omega^b 
$$
and it is of signature  $(n-2, 0)$, that is, the form 
$\widetilde{g}$ is positive semidefinite on $V^{n-1}$.

Taking exterior derivative of equation (18) by means of (2), 
we obtain the exterior quadratic equation 
\begin{equation}\label{eq:19}
\omega^a \wedge \omega_a^n =0.
\end{equation}
Applying Cartan's lemma to this equation, 
we find that 
\begin{equation}\label{eq:20}
\omega_a^n = \lambda_{ab} \omega^b, \;\; \lambda_{ab} 
= \lambda_{ba}. 
\end{equation}
The quantities $\lambda_{ab}$ are defined in a second-order 
neighborhood of a point \newline $x \in V^{n-1}$. 
It follows from equations (17) and (20) that 
\begin{equation}\label{eq:21}
\omega_1^a =  \lambda_b^a \omega^b,  
\end{equation}
where $ \lambda_b^a = g^{ac} \lambda_{cb}$, and 
 $g^{ab}$ is the inverse tensor of the tensor $g_{ab}$. 

Let us prove the existence theorem for lightlike hypersurfaces.

\begin{theorem} Lightlike hypersurfaces on a 
manifold $(M, c)$ exist, and the solution of a system 
defining such hypersurfaces depends on one function 
of $n - 2$ variables.
\end{theorem}

{\sf Proof}. The hypersurfaces in question are defined by 
equation (18) whose exterior differentiation leads to 
exterior quadratic equation (19). Equation (19) 
contains only the basis forms $\omega^a$ and does not contain the 
form $\omega^1$. This implies that for proving the existence 
we must consider only the   forms $\omega^a$ as basis 
forms of an integral manifold. The number of these 
forms is $n - 2$. If we apply to equation (19) 
the Cartan test (see [BCGGG 91] or [AG 93], pp. 12--13), 
we find that the characters of this equation are 
$s_1 = s_2 = \ldots = s_{n-2} = 1$, and the Cartan number 
is $Q = s_1 + 2 s_2 + \ldots + (n - 2) s_{n-2} 
= \frac{(n-1)(n-2)}{2}$. A general integral element 
of equation (19) depends on $N$ arbitrary parameters, where 
$N$ is the number of independent coefficients $\lambda_{ab}$ 
in equations (20). Since $\lambda_{ab} 
= \lambda_{ba}$, we have $N =  \frac{(n-1)(n-2)}{2}$. 
Thus we have $N = Q$. This proves Theorem 1. \rule{3mm}{3mm}

An {\em isotropic geodesic} on the manifold $(M, g)$ 
is a geodesic that is tangent to the isotropic cone $C_x$ at 
each of its points $x$. 

As was proved in [AG 96] (see also [AG 97]), isotropic geodesics 
are invariant with respect to a conformal transformation 
of the metric $g$. 
 We will prove now the following theorem:

\begin{theorem} A  lightlike hypersurface 
$V^{n-1} \subset M (c)$ carries a foliation formed 
 by isotropic geodesics.
\end{theorem}

{\sf Proof}. Since the straight lines $A_0 A_1$ 
are tangent to  isotropic lines on a hypersurface $V^{n-1}$, 
the equations of the isotropic foliation on $V^{n-1}$ 
have the form
\begin{equation}\label{eq:22}
\omega^a = 0.  
\end{equation}

Let us prove that the curves belonging to this foliation 
are isotropic geodesics. It is known (see [AG 96], Section 
{\bf 4.2}) that the equations of geodesics in any of 
pseudo-Riemannian metrics compatible with 
the  $CO (n-1, 1)$-structure 
can be written as 
\begin{equation}\label{eq:23}
d \omega^i + \omega^j  \omega^i_j = \alpha \omega^i, \;\;\;\;\; 
i, j = 1, \ldots , n,
\end{equation}
where $\alpha$ is a 1-form, and $d$ is the symbol of ordinary 
(not exterior) differentiation. 
In our moving frame equations (23) take the form
\begin{equation}\label{eq:24}
\renewcommand{\arraystretch}{1.3}
\left\{
\begin{array}{ll}
d \omega^1 + \omega^1  \omega_1^1 + \omega^a \omega_a^1 
+ \omega^n  \omega_n^1= \alpha \omega^1, &\\
d \omega^a + \omega^1  \omega_1^a + \omega^b \omega_b^a 
+ \omega^n  \omega_n^a = \alpha \omega^a, &
a, b = 2, \ldots , n-1, \\
d \omega^n + \omega^1  \omega_1^n + \omega^b \omega_b^n 
+ \omega^n  \omega_n^n = \alpha \omega^n. &
\end{array}
\right.
\renewcommand{\arraystretch}{1}
\end{equation}
By means of equations (18) and (21), which are valid on 
$V^{n-1}$, and equation (22) defining the isotropic 
foliation on $V^{n-1}$, the last two equations (24) are 
identically satisfied, and the remaining first equation 
determines $d\omega^1$ on a geodesic. \rule{3mm}{3mm}

Note that for lightlike hypersurfaces of 
a pseudo-Riemannian space, a similar result in a slightly 
different  terminology is given in 
[DB 96], p. 86. 

{\bf 2.} Consider the development of isotropic 
geodesics of a hypersurface $V^{n-1}$ defined by equation (22) on 
the hyperquadric $Q_1^n$. By means of (18) and (22), it follows 
from equations (12) that 
$$
d A_0 = \omega^0_0 A_0 + \omega^1_0 A_1, \;\;
d A_1 = \omega_1^0 A_0 + \omega^1_1 A_1.
$$
These equations prove that under the development, 
to the isotropic geodesics defined by equation (22) 
there corresponds an open part of the rectilinear generator 
$A_0 A_1$ of the hyperquadric $Q_1^n$. We assume that 
{\em the isotropic geodesics of the hypersurface $V^{n-1}$ 
are prolonged in such a way that they are mapped onto 
 the entire rectilinear generator $A_0 A_1$} which is 
a projective straight
 line $l$. From equation (22) it follows also that 
the family of isotropic geodesics on  a hypersurface $V^{n-1}$ 
depends on $n - 2$ parameters, and the forms $\omega^a$ 
are independent linear combinations of differentials of 
these parameters. This implies the 
following theorem:

\begin{theorem} A  lightlike hypersurface 
$V^{n-1}$ of a differential manifold $M, \dim M = n$, 
endowed with a pseudoconformal structure $CO (n-1, 1)$ is 
the image of the product $M^{n-2} \times l$, 
where $l$ is a projective straight line, under a mapping 
$f$ of this product onto the manifold $M$, 
$V^{n-1} = f (M^{n-2} \times l)$.
\end{theorem}

Note also that since the isotropic geodesics of 
a lightlike hypersurface  $V^{n-1}$ are the images 
of projective straight lines, then one can introduce 
projective coordinates both homogeneous and 
nonhomogeneous. In what follows we will use 
this remark.

Consider a displacement of the isotropic geodesic 
$l = A_0 A_1$  on a lightlike hypersurface  $V^{n-1}$. 
From equations (12), (17), and (18) it follows that 
\begin{equation}\label{eq:25}
\renewcommand{\arraystretch}{1.3}
\left\{
\begin{array}{ll}
d A_0 = \omega^0_0 A_0 + \omega^1_0 A_1 + \omega^a_0 A_a, \\
d A_1 = \omega_1^0 A_0 + \omega^1_1 A_1 + \omega_1^a A_a,
\end{array}
\right.
\renewcommand{\arraystretch}{1}
\end{equation}
where the forms $\omega_1^a$ have expressions (21). 
Consider a point $Z = A_1 - s A_0$ on the straight line 
$l$. From equations (25) and (21) it follows that 
\begin{equation}\label{eq:26}
d Z \equiv (\lambda_b^a - s \delta_b^a) \omega^b A_a 
\pmod{A_0, A_1}.
\end{equation}
The matrix $(J_b^a) = (\lambda_b^a - s \delta_b^a)$ 
is the Jacobi matrix of the mapping $f$, and its determinant,
$$
J = \det (\lambda_b^a - s \delta_b^a)
$$
is the Jacobian of this mapping. 

Since the quasiaffinor $\lambda_b^a = g^{ac} \lambda_{cb}$ 
is symmetric, its characteristic equation 
\begin{equation}\label{eq:27}
\det (\lambda_b^a - s \delta_b^a) = 0
\end{equation}
has $n - 2$ real roots if each of them is counted as many times 
as its multiplicity. This implies the following theorem.

\begin{theorem} Any isotropic geodesic $l$ of 
a lightlike hypersurface  $V^{n-1}$ of a manifold $M$ endowed 
with a pseudoconformal structure of Lorentzian signature 
carries $n - 2$ real singular points 
if each of them is counted as many times 
as its multiplicity.
\end{theorem}

{\sf Proof.} The tangent subspace to 
a lightlike hypersurface  $V^{n-1}$ at a point $Z$ 
is a subspace of the space $T_x (M)$. By (25) and (26), 
this subspace is determined by the point $Z, A_1$, and 
$C_b = (\lambda_b^a - s \delta_b^a) A_a$. If the Jacobian $J$ is 
different from 0, then these points are linearly independent and 
determine the $(n-1)$-dimensional tangent subspace 
$T_Z (V^{n-1})$. In this case a point $Z$ is a regular 
point of the hypersurface  $V^{n-1}$. If at a point 
$Z \in A_0 A_1$ the Jacobian $J$ is equal to 0, then at this 
point $\dim T_Z (V^{n-1}) < n - 1$, and this point is 
a singular point of $V^{n-1}$. The coordinates $s$ of these 
singular points can be found from equation (27) which has 
$n - 2$ real roots. \rule{3mm}{3mm}

It is obvious that the point $x = A_0$ is a regular point 
of the straight line $l$.

Denote by $s_a$ the roots of  equation (27). Then the singular 
points of the isotropic geodesic $l$ have the expressions
\begin{equation}\label{eq:28}
F_a = A_1 - s_a A_0.
\end{equation}
In the paper [AG 98b], for a lightlike hypersurface 
of a pseudo-Riemannian de Sitter space 
we investigated the structure of these singular points 
 and the structure of $V^{n-1}$ itself taking into account 
multiplicities of singular points.  Many of the results of 
[AG 98b] are still valid for a lightlike hypersurface  $V^{n-1}$ 
on a manifold endowed with a pseudoconformal structure.

One more important property 
of  a lightlike hypersurface  $V^{n-1}$ of a space 
with a pseudoconformal structure follows from our 
considerations. This property is described in the following 
theorem.

 \begin{theorem} Under the development of 
a lightlike hypersurface  $V^{n-1}$ of a manifold $M$ endowed 
with a pseudoconformal structure of Lorentzian signature onto a 
hyperquadric $Q_1^n \subset P^{n+1}$ along its  isotropic 
geodesic $l$, to the tangent hyperplanes $T_Z (V^{n-1})$ 
at regular points $Z$ of the line $l$, there corresponds a 
unique subspace $T_l$ of dimension $n - 1$ that is tangent 
to the hyperquadric $Q_1^n$ at all points of the line $l$.
\end{theorem}

{\sf Proof.} In fact, from (25) and (26) it follows that 
at  regular points $Z$ of the line $l$, i.e., for $J\neq 0$, 
the tangent subspace $T_Z$ is determined by 
the same points $A_0, A_1, A_2, \ldots , A_{n-1}$. 
Therefore this subspaces are not changed when 
a point $x$ moves along the line $l$. \rule{3mm}{3mm}

\section{The fundamental geometric objects and 
\newline tensors of a lightlike hypersurface \newline  
defined in a second-order neighborhood}

{\bf 1.} 
Singular points $F_a$  are defined invariantly on an isotropic 
geodesic $l$ of a lightlike hypersurface  $V^{n-1}$. But the 
coordinates $s_a$ of these points depend on the choice of the 
points $A_0$ and $A_1$ on this isotropic geodesic and on 
normalization of these points. The coefficients of 
characteristic equation (27) also depend on the choice and 
normalization of these two points. The point $A_0$ can move 
freely along the straight line $l$, since $A_0$ is an arbitrary 
point of a lightlike hypersurface  $V^{n-1}$. 
 A displacement of 
this point is determined by a parameter $u^1$ whose 
differential $du^1$ is contained in the basis form 
$\omega^1$. As to the point $A_1$, its freedom of motion 
can be restricted. For example, we can suppose that 
the point $A_1$ is the {\em harmonic pole} 
(introduced in [Cas 50]) of 
the point $A_0$ with respect to the foci $F_a$ of 
the isotropic geodesic $l$. Then the displacement 
of the point $A_1$ is determined by the same parameter 
$u^1$ which determines the displacement of $A_0$. 

The coordinate $\lambda$ of the harmonic pole 
of the point $A_0$ with respect to the foci $F_a$ is equal to 
the arithmetic mean of coordinates of the foci $F_a$:
$$
\lambda = \frac{1}{n-2} \sum_{a} s_a.
$$
But the sum of the roots of algebraic equation (27) is 
the negative of the coefficient in $s^{n-3}$ of 
this equation, that is, this sum is the trace of 
the quasiaffinor $\lambda_b^a$. Thus
\begin{equation}\label{eq:29}
\lambda = \frac{1}{n-2} \lambda_a^a 
= \frac{1}{n-2}  \lambda_{ab} g^{ab}, 
\end{equation}
and we have the following expression of the harmonic pole $H$:
\begin{equation}\label{eq:30}
H = A_1 - \lambda A_0 
\end{equation}
(see Figure 3). 

\vspace*{35mm}

\begin{center}
Figure 3
\end{center}

If we superpose the point $A_1$ with the harmonic pole $H$, 
the we obtain $\lambda = 0$. After such a normalization, 
all remaining coefficients of characteristic equation (27) 
become relative invariants of the  hypersurface  $V^{n-1}$. 
The weights of these invariants are equal to the degrees 
of a component 
of quasiaffinor $\lambda_b^a$ that occurs in the expressions 
of these components.

{\bf 2.} Consider the second prolongation of the basic 
differential equations (18) of a lightlike hypersurface  
$V^{n-1}$ in a manifold $M$ endowed with the pseudoconformal 
$CO (n-1, 1)$-structure of Lorentzian signature. To this end,  
using equations (3)--(5), 
we take exterior derivatives of equations (20) obtained 
in the first prolongation of equations (18).
 As a result,  we 
arrive at the following exterior quadratic equations:
\begin{equation}\label{eq:31}
[\nabla \lambda_{ab} - \lambda_{ab}(\omega_0^0 + \omega_1^1) 
- g_{ab} \omega^0_1 + (2C_{ab1}^n + \lambda_{al} g^{lc} 
\lambda_{cb}) \omega^1 + C^n_{abc} \omega^c] \wedge \omega^b = 0,
\end{equation}
where $\nabla \lambda_{ab} = d\lambda_{ab} - \lambda_{cb} 
(\omega^c_a - \delta^c_a \omega_0^0) - \lambda_{ac} 
(\omega^c_b - \delta^c_b \omega_0^0)$. 
 By Cartan's lemma, we find from (31) that
\begin{equation}\label{eq:32}
\nabla \lambda_{ab} - \lambda_{ab}(\omega_0^0 + \omega_1^1) 
- g_{ab} \omega^0_1 + (2C_{ab1}^n + \lambda_{al} g^{lc} 
\lambda_{cb}) \omega^1 + C^n_{abc} \omega^c
= \lambda_{abc} \omega^c,
\end{equation}
where $\lambda_{abc}$ are symmetric with respect to all lower 
indices. 

The quantities $C^n_{ab1}$ are symmetric with 
respect to the indices $a$ and $b$ since by (7) and (8) we have 
$$
C^n_{ab1} = - C_{1ab1} = - C_{b11a} = - C_{1ba1} = C^n_{ba1}.
$$ 
If we alternate equations (32) with respect to 
the indices $a$ and $b$, then we find that 
$C^n_{[ab]c} = 0$. This implies that $C^n_{abc} = C^n_{bac}$. 
By (7) and (8) we have 
$
C^n_{abc} = - C^n_{acb}. 
$
It follows that 
$$
C^n_{abc} = - C^n_{acb} = - C^n_{cab} = C_{cba}^n = 
 C^n_{cab} = -  C^n_{bac} = -  C^n_{abc}.
$$
Thus 
the components $C^n_{abc}$ of the curvature tensor  satisfy 
the conditions
$$
C^n_{abc} = 0.
$$
As a result, equation (32) takes the form 
\begin{equation}\label{eq:33}
\nabla \lambda_{ab} - \lambda_{ab}(\omega_0^0 + \omega_1^1) 
- g_{ab} \omega^0_1 + (2C_{ab1}^n + \lambda_{al} g^{lc} 
\lambda_{cb}) \omega^1 = \lambda_{abc} \omega^c.
\end{equation}

Note also that using the operator $\nabla$, we can write 
equations (6) and the corresponding equations for the 
tensor $g^{ab}$ in the form
\begin{equation}\label{eq:34}
\nabla g_{ab} = 2 g_{ab} \omega^0_0, \;\; 
\nabla g^{ab} = - 2 g^{ab} \omega^0_0.
\end{equation}

For a fixed point $x \in V^{n-1}$ (i.e., for 
$\omega^1 = \omega^a = 0$), we find from (32) that 
\begin{equation}\label{eq:35}
\nabla_\delta \lambda_{ab} - \lambda_{ab}(\pi_0^0 + \pi_1^1) 
- g_{ab} \pi_1^0 = 0,
\end{equation}
where $\delta = d|_{\omega^1 = \omega^a = 0}, \;
\pi^i_j = \omega^i_j (\delta) = \omega^i_j|_{\omega^1 = \omega^a = 0}$, 
and $\nabla_\delta \lambda_{ab} = \delta \lambda_{ab} 
- \lambda_{cb} (\pi^c_a - \delta^c_a \pi_0^0) 
- \lambda_{ac} (\pi^c_b - \delta^c_b \pi_0^0).$

Equations (35) prove that the quantities $\lambda_{ab}$ 
do not form a tensor since they are changed under a displacement 
of the point $A_1$ along the isotropic geodesic $l = A_0 A_1$. 
 However, the quantities $\{g_{ab}, \lambda_{ab}\}$ allow us to 
construct a tensor  defined in a 
second-order differential neighborhood of 
a point $x \in V^{n-1}$. To this end, we consider 
the geometric object $\lambda$ defined by formula (29). 
We set $\omega^1 = \omega^a = 0$ and 
differentiate (29), using (35) and the relation 
$\nabla_\delta g^{ab} = - 2 g^{ab} \pi_0^0$ 
(which follows from (34)). 
As a result, we find that 
$\lambda$ satisfies the following differential equation:
\begin{equation}\label{eq:36}
\delta \lambda  + \lambda (\pi_0^0 - \pi_1^1) - \pi_1^0 = 0.
\end{equation}
Using the quantities $\lambda_{ab}, \lambda$, and $g_{ab}$, 
we construct the quantities
\begin{equation}\label{eq:37}
h_{ab} = \lambda_{ab} - \lambda g_{ab}.
\end{equation}
Differentiating (37) with respect to the fiber parameters
 and using (35), (36) and (34), 
we find that these new quantities satisfy the following 
differential equations:
\begin{equation}\label{eq:38}
\nabla_\delta h_{ab}  =  h_{ab} (\pi_0^0 + \pi_1^1),
\end{equation}
where $\nabla_\delta h_{ab}$ has the expression similar 
to that of $\nabla_\delta \lambda_{ab}$. 
It follows from (38) that the quantities $h_{ab}$ form 
a symmetric relative (0, 2)-tensor that is defined in a 
second-order differential neighborhood of 
a point $x \in V^{n-1}$. Contracting equation (37) 
with the tensor $g^{ab}$, we find that 
\begin{equation}\label{eq:39}
 h_{ab}  g^{ab} = 0,
\end{equation}
i.e., tensor $h_{ab}$ is apolar to the tensor $g_{ab}$.

By means of the tensor $h_{ab}$, we can construct 
the affinor 
\begin{equation}\label{eq:40}
h_b^a =   g^{ac}  h_{cb} = \lambda_b^a - \lambda \delta_b^a,
\end{equation}
that is also defined in a 
second-order differential neighborhood of 
a point $x \in V^{n-1}$. This affinor is trace-free, since 
it is easy to check that $h_a^a = 0$.

{\bf 3.} Equations (37) prove that if 
we superpose the point $A_1$ with the harmonic pole $H$ 
of the point $A_0$ with respect to 
the foci $F_a$ of the isotropic geodesic $l$, then the 
quantities $\lambda_{ab}$ will be identically equal to 
the corresponding components of the tensor $h_{ab}$. 
It is naturally to call the tensor $h_{ab}$ 
the second fundamental tensor of the hypersurface 
$V^{n-1}$, and the affinor $h_b^a$ the {\em Burali--Forti 
affinor} of $V^{n-1}$ (cf. [Bu 12]). (Note that the authors 
of [DB 96] called $h_b^a$ the shape operator; see 
 [DB 96], $\;$ pp. 85, 154, and 160.)  

The following two theorems clarify a geometric meaning 
of the affinor $h_b^a$.

\begin{theorem} The harmonic pole $H$ of the point $A_0 =x$ 
with respect to the foci $F_a$ of the isotropic geodesic 
$l = A_0 A_1$ is its regular point if and only if 
 $h = \det (h_b^a) \neq 0$.  
\end{theorem}

{\sf Proof.} Superpose the vertex $A_1$ of the frame associated 
with an isotropic geodesic $l = A_0 A_1$ with the 
harmonic pole $H$ of its point $A_0$, $A_1 = H$. Then by 
(30), (37), (40),  and (39), we find that 
\begin{equation}\label{eq:41}
\lambda = 0, \;\; \lambda_{ab} = h_{ab}, \;\; 
\lambda_b^a =  h_b^a, \;\; h_a^a = 0,
\end{equation}
and by (21), we obtain  
\begin{equation}\label{eq:42}
\omega_1^a =  h_b^a \omega^b.
\end{equation}
Hence from (25) it follows that 
\begin{equation}\label{eq:43}
d A_1 = \omega_1^0 A_0 +  \omega^1_1 A_1 
+ h_a^b \omega^a A_b.
\end{equation}
The tangent subspace to the hypersurface $V^{n-1}$ 
at the point $A_1 = H$ is determined by the points 
$A_1, A_0$, and $C_a = h_a^b A_b$, and if $h \neq 0$, 
then $\dim \; T_H (V^{n-1}) = n - 1$. This implies 
 Theorem 6. 
 \rule{3mm}{3mm}

A point $x$ of a lightlike hypersurface 
$V^{n-1} \subset (M, c)$ is called {\em umbilical} if 
$h_b^a = 0$ at this point. A hypersurface 
$V^{n-1}$ is called {\em totally umbilical} if all its 
points are umbilical, i.e., if the tensor $h_b^a$ 
vanishes on $V^{n-1}$.

\begin{theorem} Every  isotropic geodesic of a 
totally umbilical hypersurface 
$V^{n-1} \subset (M, c)$ carries  a single 
$(n - 2)$-fold singular point coinciding with its 
 harmonic pole $H$. Moreover, 
\begin{description}
\item[(a)] If a manifold $(M, c)$ is not conformally 
flat, then the point $H$ describes a singular curve 
on $V^{n-1}$.

\item[(b)] If a manifold $(M, c)$ is  conformally 
flat, then the point $H$ is fixed. In this case 
the image of the hypersurface $V^{n-1}$ 
on a hyperquadric $Q^n_1$ is the isotropic 
cone $C_H = Q^n_1 \cap T_H (Q^n_1).$

\item[(c)] Conversely, any  isotropic 
cone $C_y$ of the conformal space $C^n_1$ is a 
totally umbilical lightlike hypersurface.
\end{description}
\end{theorem}

{\sf Proof.} First note that for $h_b^a = 0$, it follows from 
(42) that 
$$
\lambda_b^a = \lambda \delta_b^a.
$$
This allows us to write equation (27) defining 
the coordinates of singular points of the line $l$ 
in the form
$$
(s - \lambda)^{n-2} = 0.
$$
It follows that the harmonic pole $H = A_1 - \lambda A_0$ 
coincides with a single $(n-2)$-fold singular point 
of the line $l$ and is not changed when the point $A_0$ moves 
along the line $l$. 

Next, on  a totally umbilical lightlike hypersurface $V^{n-1}$ 
equations (42) take the form 
\begin{equation}\label{eq:44}
 \omega_1^a = 0.
\end{equation}
Taking exterior derivatives of equations (44) and using 
 (4), we find that 
\begin{equation}\label{eq:45}
 \omega_1^0 \wedge  \omega_0^a + 2 C^a_{11b} 
\omega^1 \wedge  \omega^b + C^a_{1bc} \omega^b 
\wedge  \omega^c= 0.
\end{equation}
From (45) it follows that 1-form $\omega_1^0$ is 
a linear combination of the basis forms $\omega^1$ 
and $\omega^a$:
$$
\omega_1^0 = \mu \omega^1 + \mu_a \omega^a.
$$
Using equations (45), it is easy to prove 
that
$$
\mu = 0, \;\; \mu_a = \frac{2}{n-3} C^b_{1ba}.
$$
Of course, we should assume that  $n \geq 4$.

By (44), we have 
$$
d A_1 = \omega_1^0 A_0 +  \omega^1_1 A_1.
$$
It follows that the point $A_1$ describes a curve tangent to 
the isotropic geodesic $l = A_0 A_1$. This proves 
part (a) of Theorem 7.

If a manifold $(M, c)$ is conformally flat, then it 
follows from (45) that $\omega_1^0 = 0$ and 
for $n \geq 4$,  equation (44) takes the form
\begin{equation}\label{eq:46}
d A_1 =  \omega^1_1 A_1,
\end{equation} 
This implies that the point $A_1 = H$ is fixed. 
Under the Darboux mapping 
on a hyperquadric $Q^n_1$ of the space $P^{n+1}$,  
the isotropic geodesics of the  hypersurface  $V^{n-1}$ 
are mapped into rectilinear generators of the isotropic 
cone $C_y$ whose vertex is the image of the point $A_1$ 
under the Darboux mapping, $y = A_1$. 
This proves part (b) of Theorem 7. 

Finally, we will prove part (c) of Theorem 7. Let 
$x$ be an arbitrary point of an isotropic cone $C_y$ 
with vertex $y$. With this cone we associate 
a first-order frame bundle in such a way that $A_0 = x$ 
and $A_1 = y$. Then since the point $y$ is fixed, 
we have equations (46). They imply $\omega_1^a = 0, \;
\lambda_b^a = 0, \; h_b^a = 0$. This proves part (c).
\rule{3mm}{3mm}

\section{An affine connection on a lightlike \newline hypersurface}

Let us find conditions 
under which on a lightlike hypersurface  $V^{n-1}$ of a 
 manifold endowed with a pseudoconformal structure  
of Lorentzian signature there will be defined an affine 
connection. Such a hypersurface is defined by equations 
(17) and (18), and its basis forms are $\omega^1, \; 
\omega^a, \; a = 2, \ldots , n-1$. Therefore, 
on such a hypersurface equations (2) take the form
\begin{equation}\label{eq:47}
\renewcommand{\arraystretch}{1.3}
\left\{
\begin{array}{ll}
d \omega^1 = \omega^1 \wedge (\omega_1^1 - \omega^0_0) 
+ \omega^a \wedge \omega^1_a, \\
d \omega^a = \omega^1 \wedge \omega_1^a + \omega^b 
\wedge (\omega_b^a  - \delta_b^a \omega^0_0).
\end{array} 
\right.
\renewcommand{\arraystretch}{1}
\end{equation}
Thus the matrix 1-form
$$
\omega = \pmatrix{\omega_1^1 - \omega^0_0 &  \omega_a^1 \cr
                 \omega_1^a & \omega_b^a  - \delta_b^a\omega^0_0} 
$$
defines on $V^{n-1}$ an affine structure. To define an 
 affine connection, the form $\omega$ must satisfy the 
structure equation
$$
d\omega +  \omega \wedge \omega = \Omega,
$$
where $\Omega$ is the curvature 2-form of this connection 
which is a linear combination of exterior products 
of the basis forms $\omega^1$ and $\omega^a$ 
(see, for example, [KN 63], Ch. III).

We  differentiate the form $\omega$ componentwise and 
apply equations (3) and (4). 
As a result, we find that 
\begin{equation}\label{eq:48}
\renewcommand{\arraystretch}{1.3}
\left\{
\begin{array}{ll}
d (\omega_1^1 - \omega^0_0) -  \omega_1^a \wedge \omega^1_a 
= 2 \omega_1^0 \wedge \omega_0^1 +\omega^a \wedge \omega_a^0  
+ C^1_{1kl} \omega^k \wedge \omega^l, \\
d \omega_a^1 - \omega^1_a \wedge (\omega^1_1 -\omega_0^0) 
-  (\omega_a^b - \delta_a^b \omega^0_0) \wedge \omega_b^1  \\
= \omega_a^0 \wedge \omega_1^0  - g_{ab} \omega^b \wedge \omega_n^0  
+ C^1_{akl} \omega^k \wedge \omega^l, \\
d \omega_1^a - (\omega^1_1 -\omega_0^0) \wedge \omega_1^a 
- \omega_1^b \wedge (\omega_b^a - \delta_b^a \omega^0_0)  \\
= \omega_1^0 \wedge \omega_0^a   
+ C^a_{1kl} \omega^k \wedge \omega^l, \\
d (\omega_b^a - \delta_b^a \omega^0_0) -  \omega_b^1 
\wedge \omega^a_1 - (\omega_b^c - \delta_b^c \omega_0^0) 
\wedge (\omega_c^a  - \delta_c^a \omega_0^0) \\ 
= \omega^0_b \wedge \omega^a_0  + \omega_b^n \wedge \omega_n^a 
+ g_{bc} g^{ae} \omega^c \wedge \omega_e^0   
- \delta_b^a(\omega_0^1 \wedge \omega_1^0 
+ \omega_0^c \wedge \omega_c^0)  
+ C^a_{bkl} \omega^k \wedge \omega^l; 
\end{array} 
\right.
\renewcommand{\arraystretch}{1}
\end{equation}
in these formulas $a, b = 2, \ldots , n - 1; \, 
k, l = 1, 2, \ldots , n - 1$. 

The right-hand sides of these equations are not expressed yet 
in terms of basis forms since the 1-forms $\omega^0_1, \;
\omega^0_a, \; \omega^0_n$, and $\omega^a_n = g^{ab} \omega^1_b$ 
are fiber forms in the first-order frame bundle associated with 
a lightlike hypersurface  $V^{n-1}$. To make these forms 
principal, it is necessary to specialize our moving frames. 
The forms $\omega_a^1$ become principal if on $V^{n-1}$ 
there is given a screen distribution $S$, and along 
with distribution, also a field of normalizing isotropic lines 
$A_0 A_n$ is given (see Section {\bf 2.1}). As it was indicated 
in [DB 96], such a specialization is sufficient for 
defining an affine connection on a lightlike hypersurface  
$V^{n-1}$ of a pseudo-Riemannian manifold $(M, g)$ of Lorentzian 
signature. 

However, in order to define an affine connection 
on a manifold $(M, c)$ endowed with a pseudoconformal 
structure of Lorentzian signature, it is not 
sufficient to have only a screen distribution. 
For the 1-forms $\omega^0_1, \; \omega^0_a$, 
and $\omega^0_n$ to be principal, it is necessary that 
in the tangent space $T_x (M)$ of a point $x \in V^{n-1}$, 
a hyperplane $L_x$ not passing through the point $x$ 
is given. Then a field $L$ of such hyperplanes $L_x$ 
together with 
a screen distribution define an affine connection on $M$. 
Thus the following theorem is valid.

 \begin{theorem} To define an affine connection 
on a lightlike hypersurface  $V^{n-1}$  
of a manifold $(M, c)$,  it is  
sufficient to assign on $V^{n-1}$ a screen distribution $S$ and a 
field $L$ of normalizing hyperplanes $L_x$ belonging to 
the tangent bundle $T (M)$.
\end{theorem}

Note that on a pseudo-Riemannian manifold $(M, g)$ 
the role of normalizing hyperplanes is played by 
the planes at infinity of the tangent spaces $T_x (M)$. 

In what follows, on an isotropic hypersurface $V^{n-1}$,
 we will make an invariant construction of a screen distribution 
$S$ and a field  $L$ of normalizing hyperplanes that 
are intrinsically connected with the geometry of 
the  hypersurface $V^{n-1}$.

\section{Construction of the main part of an 
invariant normalization of a lightlike hypersurface}

{\bf 1.} We will derive now some formulas that 
will be used later. Differentiating equation (29) and 
applying equation (31), we obtain the following 
Pfaffian equation:
\begin{equation}\label{eq:49}
\renewcommand{\arraystretch}{1.3}
\begin{array}{ll}
d \lambda + \lambda (\omega_0^0 - \omega^1_1) -  \omega_1^0 
= & - \displaystyle 
\frac{1}{n-2} g^{ab} (\lambda_{ae} g^{ec} \lambda_{cb} 
+ 2 C^n_{ab1} ) \wedge \omega^1 \\
&+  \displaystyle 
\frac{1}{n-2} g^{ab} \lambda_{abc}  \wedge \omega^c.\end{array} 
\renewcommand{\arraystretch}{1}
\end{equation}
From the last equation of (8) it follows that 
$$
g^{ab} C^n_{ab1} = 0, \;\; C_{n11c} = 0.
$$
This allows us to write equation (49) in the form:
\begin{equation}\label{eq:50}
d \lambda + \lambda (\omega_0^0 - \omega^1_1) -  \omega_1^0 
=- \displaystyle \frac{1}{n-2} g^{ab} \lambda_{ae} g^{ec} 
\lambda_{cb} \omega^1 
+ \displaystyle  \frac{1}{n-2} g^{ab} \lambda_{abc}  \omega^c.
\end{equation}
Set 
\begin{equation}\label{eq:51}
\mu = \frac{1}{n-2} g^{ab} \lambda_{ae} g^{ec} \lambda_{cb}, \;\;\;\; 
\mu_c = - \frac{1}{n-2}  g^{ab} \lambda_{abc}.
\end{equation}
Then equation (50) can be written as 
\begin{equation}\label{eq:52}
d \lambda + \lambda (\omega_0^0 - \omega^1_1) -  \omega_1^0 
= - \mu \omega^1 - \mu_a \omega^a.
\end{equation}
Note that {\em the quantities $\mu$ and $\mu_a$ are 
 defined in a second- and third-order neighborhood, respectively, 
 of a point $x$ of the hypersurface $V^{n-1}$}.

Note also that since the tensor $g_{ab}$ is positive 
definite, the quantity $\mu$ defined by (50) is equal to 0 at a 
point $x$ if and only if the quasitensor $\lambda_{ab}$ is equal 
to 0. But then equations (37) imply that the tensor $h_{ab}$ 
is equal to 0 at the point $x$, and as a result, the point $x$ 
is umbilical. Thus {\em the quantity $\mu$ is equal to $0$ 
at umbilical points of the hypersurface $V^{n-1}$, and only at 
such points.} In what follows we will assume that the 
hypersurface 
$V^{n-1}$ does not have umbilical points and that $\mu \neq 0$ 
on $V^{n-1}$. Moreover, from (51) it follows that 
$\mu > 0$.

If we take exterior derivative of equation (52), we 
obtain the following exterior quadratic equation:
\begin{equation}\label{eq:53}
\renewcommand{\arraystretch}{1.3}
\begin{array}{ll}
&[d \mu + 2 \mu (\omega_0^0 - \omega^1_1) - 2\lambda  \omega_1^0] 
\wedge \omega^1_0 \\
+& [\nabla \mu_a +  \mu_a (\omega_0^0 - \omega^1_1) + 
h_a^b \omega_b^0 - (\mu \delta_a^b - \lambda \lambda_a^b)  
\omega^1_b \\
+ & \mu_b \lambda_a^b \omega_0^1 
-  2(\lambda C_{11a}^1 + C_{11a}) \omega_0^1 
+ (\lambda C^1_{1ab} + C_{1ab}) \omega^b] \wedge \omega^a = 0,
\end{array} 
\renewcommand{\arraystretch}{1}
\end{equation}
where $\nabla \mu_a = d \mu_a - \mu_b (\omega_a^b - \delta_a^b 
\omega_0^0)$.  
 Applying Cartan's lemma to equation (53), we find 
that
\begin{equation}\label{eq:54}
\renewcommand{\arraystretch}{1.3}
\left\{
\begin{array}{ll}
d \mu \!\!\!\! &+ 2 \mu (\omega_0^0 - \omega^1_1) 
- 2 \lambda  \omega_1^0 
= \nu \omega^1 + \nu_a  \omega^a, \\
\nabla \mu_a \!\!\!\!&+  \mu_a (\omega_0^0 - \omega^1_1) + 
h_a^b \omega_b^0 - (\mu \delta_a^b - \lambda \lambda_a^b)  
\omega^1_b +  \mu_b \lambda_a^b \omega_0^1 \\
\!\!\!\!&- 2 (\lambda C_{11a}^1 + C_{11a}) \omega_0^1 
+ (\lambda C^1_{1ab} + C_{1ab}) \omega^b 
= \nu_a \omega^1 + \nu_{ab} \omega^b,
\end{array} 
\right.
\renewcommand{\arraystretch}{1}
\end{equation}
where $\nu_{ab} = \nu_{ba}$. Here the quantities $\nu$ and 
$\nu_a$ 
are  defined in a third-order neighborhood of 
a point $x \in V^{n-1}$, and the quantities $\nu_{ab}$ 
are defined in a fourth-order neighborhood of 
 $x \in V^{n-1}$.

{\bf 2.} In what follows we assume that the point $A_1$ of an 
isotropic geodesic $l = A_0 A_1$ is superposed with the harmonic 
pole $H$ of the point $A_0$ with respect to the singular points 
$F_a$ of this straight line. Then equations (41) and (42) hold, 
equation (52) takes the form
\begin{equation}\label{eq:55}
  \omega_1^0 = \mu \omega^1 + \mu_a \omega^a,
\end{equation}
and equations (54) become 
\begin{equation}\label{eq:56}
\renewcommand{\arraystretch}{1.3}
\left\{
\begin{array}{ll}
d \mu \!\!\!\! &+ 2 \mu (\omega_0^0 - \omega^1_1) 
 = \nu \omega^1 + \nu_a \omega^a, \\
\nabla \mu_a \!\!\!\! &+  \mu_a (\omega_0^0 - \omega^1_1) 
+ h_a^b \omega_b^0  - \mu \omega_a^1 
- 2  C_{11a} \omega^1 + C_{1ab}  
\omega^b = \nu_a \omega^1 + \nu_{ab} \omega^b.
\end{array} 
\right.
\renewcommand{\arraystretch}{1}
\end{equation}

Let us write the last equations for fixed 
principal parameters, i.e., for $\omega^1 = \omega^a = 0$:
\begin{equation}\label{eq:57}
\renewcommand{\arraystretch}{1.3}
\left\{
\begin{array}{ll}
\delta \mu = 2 \mu (\pi_1^1 - \pi^0_0), \\
\nabla_\delta \mu_a +  \mu_a (\pi_0^0 - \pi^1_1) 
+ h_a^b \pi_b^0  - \mu \pi_a^1= 0.
\end{array} 
\right.
\renewcommand{\arraystretch}{1}
\end{equation}
The first equation of (57) proves that after the specialization 
of moving frames described above, the quantity $\mu$ becomes a 
relative invariant of weight 2. As we showed earlier, 
this relative invariant is different from 0 at non-umbilical 
points of the hypersurface $V^{n-1}$.

By (51) the invariant $\mu$ can be written 
now in the following form:
\begin{equation}\label{eq:58}
\mu =  \frac{1}{n-2} g^{ab} h_{ae} g^{ec} h_{cb} 
 = \frac{1}{n-2} h^b_a h^a_b.
\end{equation}

The second equation of (57) contains two groups of fiber forms,  
$\omega_a^0$ and $\omega_a^1$, and this is the reason that 
 the reduction 
of the object $\mu_a$ to 0 does not make these forms principal. 
Hence  we should also consider the object 
$\nu_a$ occurring in equations (56) which is also defined 
in a third-order neighborhood of a point $x \in V^{n-1}$. 
In order to find differential equations which this object 
satisfies, we take the exterior derivative of the first 
equation of (56). As a result, we find that 
\begin{equation}\label{eq:59}
\renewcommand{\arraystretch}{1.3}
\begin{array}{ll}
&\!\!\!\! 
[d \nu + 3 \nu (\omega_0^0 - \omega^1_1) - \nu_a \omega_1^a 
+ 4 \mu \omega_1^0] \wedge \omega^1_0 \\
+&\!\!\!\! [\nabla \nu_a +  2 \nu_a (\omega_0^0 - \omega^1_1) + 
2 \mu  \omega_a^0 - (2 \mu h_a^b + \nu \delta_a^b)\omega_b^1]  \wedge 
\omega^a \\
+ &\!\!\!\! 2\mu (2 C_{11a}^1 \omega^1 \wedge \omega^a +  C^1_{1ab} 
 \omega^a \wedge \omega^b) = 0.
\end{array} 
\renewcommand{\arraystretch}{1}
\end{equation}
Substituting the forms $\omega_1^a$ and $\omega_1^0$ 
in equation (59) by their values taken from equations (42) and 
(55), we obtain
\begin{equation}\label{eq:60}
\renewcommand{\arraystretch}{1.3}
\begin{array}{ll}
[d \nu \!\!\!\!&+ 3 \nu (\omega_0^0 - \omega^1_1)] \wedge \omega^1 
+ [\nabla \nu_a + 2 \nu_a (\omega_0^0 - \omega^1_1) 
+ 2\mu \omega_a^0 - (2 \mu h_a^b + \nu \delta_a^b)\omega_b^1  \\
\!\!\!\!& +  (\nu_b h^b_a - 4 \mu \mu_a) \omega^1 
+ 2\mu (2 C_{11a}^1 \omega^1  -  C^1_{1ab}  \omega^b)] \wedge \omega^a = 0,
\end{array} 
\renewcommand{\arraystretch}{1}
\end{equation}
where $\nabla \nu_a = d \nu_a - \nu_b (\omega_a^b - \delta_a^b 
\omega_0^0)$. 
Applying Cartan's lemma to equation (60), we find that 
\begin{equation}\label{eq:61}
\renewcommand{\arraystretch}{1.3}
\left\{
\begin{array}{ll}
d \nu \!\!\!\!&+ 3 \nu (\omega_0^0 - \omega^1_1) = \rho \omega^1 
+ \rho_a \omega^a, \\
\nabla \nu_a \!\!\!\!&+  2\nu_a (\omega_0^0 - \omega^1_1) + 
2 \mu  \omega_a^0 - (2 \mu h_a^b + \nu \delta_a^b)\omega_b^1 \\ 
\!\!\!\!&+  (\nu_b h^b_a - 4 \mu \mu_a) \omega^1 
+ 2\mu (2 C_{11a}^1 \omega^1  -  C^1_{1ab}  \omega^b) 
 = \rho_a \omega^1 + \rho_{ab} \omega^b.
\end{array} 
\right.
\renewcommand{\arraystretch}{1}
\end{equation}
The coefficients $\rho, \rho_a$, and $\rho_{ab}$ 
in equations (61) are connected with a fourth-order differential 
neighborhood of a point $x \in V^{n-1}$ and  
$\rho_{ab} = \rho_{ba}$.

For fixed principal parameters (i.e., for $\omega^1 = \omega^a = 0$),  
equations (61) take the form 
\begin{equation}\label{eq:62}
\renewcommand{\arraystretch}{1.3}
\left\{
\begin{array}{ll}
\delta \nu + 3 \nu (\pi_0^0 - \pi_1^1) = 0, \\
\nabla_\delta \nu_a + 2 \nu_a (\pi_0^0 - \pi_1^1) 
+ 2\mu \pi_a^0 - (2\mu h_a^b + \nu \delta_a^b) \pi_b^1 = 0.
\end{array} 
\right.
\renewcommand{\arraystretch}{1}
\end{equation}
When we derived (62), we took into account that 
by (41) and (55), $\pi_1^a = \pi_1^0 = 0$.

{\bf 3.} Using the geometric objects $\mu_a$ and $\nu_a$, 
we will construct now an invariant normalization of a lightlike 
hypersurface $V^{n-1}$.

In the tangent space $T_x (V^{n-1})$, we consider the 
subspace $R_x$ which is complementary to the straight line 
$A_0 A_1$, and take the points 
\begin{equation}\label{eq:63}
C_a = A_a + y_a A_0 + z_a A_1
\end{equation}
as basis points of this space. The subspace $R_x$ is invariant if 
and only if 
$$
\delta C_a = \sigma_a^b C_b.
$$
Differentiating equation (63) with respect to fiber 
parameters, we find that 
$$
\delta C_a = (\nabla_\delta y_a +  \pi_a^0) A_0 
+ (\nabla_\delta z_a + z_a (\pi_1^1 -\pi_0^0) + \pi_a^1) A_1   
+ \pi_a^b C_b.
$$
Thus the conditions for the subspace $R_x$ 
to  be invariant  is
\begin{equation}\label{eq:64}
\renewcommand{\arraystretch}{1.3}
\left\{
\begin{array}{ll}
\nabla_\delta y_a +  \pi_a^0 = 0, \\
\nabla_\delta z_a + z_a (\pi_1^1 -\pi_0^0) + \pi_a^1 = 0.
\end{array} 
\right.
\renewcommand{\arraystretch}{1}
\end{equation}

Next, using the geometric objects $\mu_a$ 
and $\nu_a$ defined earlier, we should construct 
normalizing geometric objects satisfying equations 
(64). Let us write one more time the equations which the 
objects $\mu_a$ and $\nu_a$ satisfy:
\begin{equation}\label{eq:65}
\renewcommand{\arraystretch}{1.3}
\left\{
\begin{array}{ll}
\nabla_\delta \mu_a +  \mu_a (\pi_0^0 - \pi^1_1) 
+ h_a^b \pi_b^0  - \mu \pi_a^1 = 0, \\
\nabla_\delta \nu_a + 2 \nu_a (\pi_0^0 - \pi_1^1) 
+ 2\mu \pi_a^0 - (2\mu h_a^b + \nu \delta_a^b) \pi_b^1 = 0.
\end{array} 
\right.
\renewcommand{\arraystretch}{1}
\end{equation}
We will try to solve these equations for 1-forms 
$\pi_a^0$ and $\pi_a^1$.
Construct the objects
\begin{equation}\label{eq:66}
M_a = h_a^b \mu_b + \frac{\nu}{2\mu} \mu_a - \frac{1}{2} \nu_a,  \;\; 
N_a =  \frac{1}{2} h_a^b \nu_b - \mu \mu_a.
\end{equation}
Differentiating these equations with respect to 
fiber parameters and applying formulas (65), (57), and (38), 
we find that 
\begin{equation}\label{eq:67}
\renewcommand{\arraystretch}{1.3}
\left\{
\begin{array}{ll}
\nabla_\delta M_a + 2 M_a (\pi_0^0 -  \pi_1^1) 
+ H_a^b  \pi_b^0 = 0, \\
\nabla_\delta N_a + 3 N_a (\pi_0^0 - \pi_1^1) 
- \mu H_a^b \pi_b^1 = 0,
\end{array} 
\right.
\renewcommand{\arraystretch}{1}
\end{equation}
where 
\begin{equation}\label{eq:68}
H_a^b = h_a^c h_c^b  + \frac{\nu}{2\mu} h_a^b - \mu \delta_a^b.
\end{equation}

We will establish some properties of the 
tensor $H_a^b$.

\begin{description}
\item[1)] {\em The tensor $H_a^b$ is of weight $2$.} In fact, 
differentiating  equation (68) with respect to 
fiber parameters, we find that 
\begin{equation}\label{eq:69}
\nabla_\delta H_a^b = 2 H_a^b (\pi_1^1 - \pi_0^0).
\end{equation}
Since equation (68) contains a relative invariant $\nu$ 
defined in a third-order neighborhood, the tensor $H_a^b$ 
is also connected with this neighborhood.

\item[2)] {\em The 
tensor $H_a^b$ is trace-free}. In fact, we have 
$$
H_a^a = h_a^c h_c^a - (n - 2) \mu,
$$
since $h_a^b$ is a trace-free tensor. But by formula (58) 
defining the invariant $\mu$, it follows that $H_a^a = 0$.

\item[3)] {\em The tensor $H_a^b$ can be reduced to a diagonal 
form simultaneously with the tensors $g_{ab}$ and $h_a^b$.} 
 In fact, since the tensor $g_{ab}$ is positive definite, and 
by (40) the tensor $h_a^b$ satisfies the condition $g_{ac} h_b^c 
= g_{bc} h_a^c$, it follows that  the tensors $g_{ab}$ and 
$h_a^b$ can be reduced simultaneously to  diagonal forms:
$$
g_{ab} = \delta_{ab}, \;\; h_a^b = h_a \delta_{ab}, \;\; 
a, b = 2, \ldots, n - 1.
$$
But now it follows from (68) that the tensor $H_a^b$ is also 
reduced to a diagonal form and has the following eigenvalues:
\begin{equation}\label{eq:70}
H_a = h_a^2 + \frac{\nu}{2\mu} h_a - \mu.
\end{equation}

\item[4)] {\em If the  tensor $H_a^b$ is  degenerate, 
then the relative invariants $\mu$ and $\nu$ of a lightlike 
hypersurface $V^{n-1}$ are connected by an algebraic 
equation, and the hypersurface is of a special type.} In fact, 
suppose that $\det (H_a^b) = 0$. Then at least one of 
the eigenvalues $H_a$ of the tensor $H_a^b$ is equal to 0. 
Suppose that $H_2 = 0$. This and  (70) imply that
$$
h_2^2 + \frac{\nu}{2\mu} h_2 - \mu = 0.
$$
Moreover since $h_2 \neq 0$, we have $\mu \neq 0$. 
The above written relation is an algebraic equation 
indicated earlier.

\item[5)] {\em If $n = 4$, then the  tensors $H_a^b$ and 
$h_a^b$ are proportional}:
\begin{equation}\label{eq:71}
H_a^b =  \frac{\nu}{2\mu} h_a^b.
\end{equation}
In fact, if $n = 4$, we have $a, b = 2, 3$, and
$$
\renewcommand{\arraystretch}{1.3}
\begin{array}{ll}
\mu = -\frac{1}{2} (h_2^2 + h_3^2), \\
H_2 = \frac{1}{2} (h_2^2 - h_3^2) + \frac{\nu}{2\mu} h_2, \;\;
H_3 = \frac{1}{2} (h_3^2 - h_2^2) + \frac{\nu}{2\mu} h_3.
\end{array} 
\renewcommand{\arraystretch}{1}
$$
Since the tensor $h_b^a$ is trace-free, we have $h_2 + h_3 = 0$, 
and the first terms in the expressions for $H_2$ and $H_3$ 
vanish. This implies equation (71). 

Formula (71) implies that for $n = 4$ the tensor $H_a^b$ is 
the zero-tensor either if the tensor $h_a^b$ is 
the zero-tensor (i.e., the hypersurface $V^3$ is umbilical) or 
if $\nu = 0$. 
\end{description}

{\bf 4.} Suppose now that the tensor $h_a^b$ is nondegenerate. 
As we proved in Theorem 7, this means that the harmonic pole 
$H = A_1$ of the point $x = A_0$ is a nonsingular point. 
The invariant $\mu$ defined by the formula (51) is 
different from 0. Suppose further that the relative tensor $H_a^b$ 
defined by formula (68) is also nondegenerate 
and denote by $\widetilde{H}_a^b$ the inverse tensor of the 
tensor $H_a^b$. By (69), the tensor $\widetilde{H}_a^b$ satisfies 
the equations
\begin{equation}\label{eq:72}
\nabla_\delta \widetilde{H}_a^b =  2 \widetilde{H}_a^b 
(\pi_0^0 - \pi_1^1).
\end{equation}

We construct also two more objects 
\begin{equation}\label{eq:73}
P_a = \widetilde{H}_a^b M_b, \;\; 
Q_a = \frac{1}{\mu}  \widetilde{H}_a^b N_b.
\end{equation}
By (67) and (69), these objects satisfy the equations
\begin{equation}\label{eq:74}
\nabla_\delta P_a + \pi_a^0 = 0, \;\; 
\nabla_\delta Q_a + Q_a (\pi_1^1 - \pi_0^0) + \pi_a^1 = 0.
\end{equation}
Comparing  (74) and (64), we see that equations 
(69) are satisfied if we set
$$
y_a = P_a, \;\; z_a = Q_a.
$$
This implies that in the tangents subspace $T_x (V^{n-1})$ 
the points
\begin{equation}\label{eq:75}
C_a = A_a + P_a A_0 + Q_a A_1
\end{equation}
define  an 
invariant $(n-3)$-dimensional subspace $R_x$ that is 
intrinsically connected with the geometry of the hypersurface 
$V^{n-1}$. The formulas (66) and (73) show that the geometric 
objects $P_a$ and $Q_a$ are expressed algebraically in 
terms of the objects $\mu_a$ and $\nu_a$ 
defined in a third-order neighborhood of a point $x \in V^{n-1}$. 
Hence the invariant subspace $R_x$ is defined in a third-order 
neighborhood of a point $x \in V^{n-1}$ too.

The invariant subspace $R_x$ defines an invariant {\em screen 
subspace} $S_x = [x, R_x]$ and a {\em complementary 
screen subspace} 
$\widetilde{S}_x = [H (x), R_x]$ which is the span 
of the harmonic pole $H (x)$ of a point $x \in V^{n-1}$ and the 
subspace $R_x$  (see Figure 4). 
On a lightlike hypersurface $V^{n-1}$ 
these subspaces define a {\em screen distribution} 
$S = \cup_{x \in V^{n-1}} S_x$ and a {\em complementary screen 
distribution} $\widetilde{S} = \cup_{x \in V^{n-1}} \widetilde{S}_x$ that 
are intrinsically connected with the geometry of 
the hypersurface $V^{n-1}$ and 
are defined in its third-order differential neighborhood. 
Note that in the above 
formulas, $x \in V^{n-1}$ means that all 
regular points $x$ are taken for which the 
harmonic points are regular too.

\vspace*{75mm}

\begin{center}
Figure 4
\end{center}

Thus we have proved the following result.

\begin{theorem} 
On a lightlike hypersurface $V^{n-1}$ of a manifold 
endowed with a pseudoconformal structure of Lorentzian signature, 
 an invariant  screen distribution 
$S$ and an invariant  complementary screen 
distribution $\widetilde{S}$ 
that are intrinsically connected with the geometry of 
 $V^{n-1}$ are defined by elements of a third-order differential 
neighborhood of a point $x \in V^{n-1}$ and can be constructed in 
the way indicated above.
\end{theorem}

Note that the problem of construction of an invariant 
normalization as well as of an affine connection for 
a lightlike hypersurface $V^{n-1}$ of a pseudo-Riemannian 
manifold $(M, g)$ of Lorentzian signature (for definition 
see [ON 83]) was 
considered by many authors (see [DB 96], Ch. 4). 
However, as far as we know, an invariant 
normalization and an affine connection intrinsically 
connected with the geometry of $V^{n-1}$  
were not considered in these papers. In [DB 96] (see pp. 
115--117) the authors consider a canonical normalization 
(canonical screen distribution) that is not invariant with 
respect to the Lorentzian transformations of 
the tangent space $T_x (M)$. A similar normalization 
was considered in [Bo 72].

\section{A congruence of normalizing straight lines}

{\bf 1}. In order to complete the construction of 
an invariant normalization of 
 a lightlike hypersurface $V^{n-1}$, we need only to construct an 
 invariant point on the 
isotropic normalizing straight line $x C_n$ 
that is conjugate to the screen subspace $S_x$. To simplify 
our construction, we superpose the vertices $A_a$ of our moving 
frame with the basis points $C_a$ of the invariant subspace $R_x$ 
that are defined by formulas (75). As a result, we find that 
$P_a = Q_a = 0$. Then equations (73) imply that $M_b = N_b = 0$. 
Since the tensor $H_a^b$ is nondegenerate, equations 
(66) imply that $\mu_a = \nu_a = 0$.

As a result of the specialization of moving frames 
we have made, the second equations of (56)  and (61) 
take the form:
\begin{equation}\label{eq:76}
\renewcommand{\arraystretch}{1.3}
\begin{array}{ll}
 h_a^b \omega_b^0 - \mu \omega_a^1 
- 2  C_{11a} \omega_0^1 + C_{1ab}  
\omega^b =  \nu_{ab} \omega^b, \\
2 \mu  \omega_a^0  - (2 \mu h_a^b + \nu \delta_a^b) \omega_b^1 
+ 2\mu (2 C_{11a}^1 \omega^1  -  C^1_{1ab}  \omega^b)  
= \rho_a \omega^1 + \rho_{ab} \omega^b,
\end{array} 
\renewcommand{\arraystretch}{1}
\end{equation}
The coefficients in the basis forms $\omega^1$ and $\omega^a$ in 
the right-hand sides of 
equations (76) are defined by elements of a fourth-order 
neighborhood of a point $x \in V^{n-1}$. 

Since $\det (H_a^b) \neq 0$, one can solve equations (76) 
with respect to the 1-forms $\omega_a^0$ and $\omega_a^1$. 
We will write these solutions in the form
\begin{equation}\label{eq:77}
\renewcommand{\arraystretch}{1.3}
\left\{
\begin{array}{ll}
\omega_a^0  = \sigma_a \omega^1 + \sigma_{ab} \omega^b, \\
\omega_a^1 =  \tau_a \omega^1 + \tau_{ab}  \omega^b.
\end{array} 
\right.
\renewcommand{\arraystretch}{1}
\end{equation}
The coefficients of these decompositions are expressed 
algebraically in terms of the coefficients of equation (76). 
Thus they are also defined by elements of a fourth-order 
neighborhood of a point $x \in V^{n-1}$. 

Taking exterior derivative of the second equation of (77), 
applying Cartan's lemma to exterior quadratic equation obtained, 
and fixing the principal parameters, we find that 
\begin{equation}\label{eq:78}
\renewcommand{\arraystretch}{1.3}
\left\{
\begin{array}{ll}
\nabla_\delta \tau_a = 0, \\
\nabla_\delta \tau_{ab} + \tau_{ab} (\pi_1^1 - \pi^0_0) 
- g_{ab} \pi_n^0 = 0. 
\end{array} 
\right.
\renewcommand{\arraystretch}{1}
\end{equation}
The first equation of (78) proves that the quantities $\tau_a$ 
form a  covector. The second equation of (78) allows us 
to construct a geometric object that fixes a point on a
 rectilinear generator $A_0 A_1$. This can be 
done as follows.

Consider the quantity
\begin{equation}\label{eq:79}
\tau = \frac{1}{n-2} \tau_{ab} g^{ab}.
\end{equation}
By means of 
the second equation of 
(78) and equation (17), 
this quantity satisfies the equation
\begin{equation}\label{eq:80}
\delta \tau + \tau (\pi_0^0 + \pi^1_1) - \pi_n^0 = 0. 
\end{equation}

Consider a point $Z = A_n + z A_0$ on the normalizing straight 
line $A_0 A_1$. Differentiating this point with respect to 
fiber parameters, we find that 
$$
\delta Z = [\delta z + z (\pi_0^0 + \pi^1_1) + \pi_n^0] A_0 
- \pi_1^1 Z. 
$$

It follows that the point $Z$ is invariant if and only if 
its coordinate $z$ satisfies the equation
\begin{equation}\label{eq:81}
\delta z + z (\pi_0^0 + \pi^1_1) + \pi_n^0 = 0.
\end{equation}
Comparing equations (80) and (81), we see that $z = - \tau$ 
is a solution of equation (81). Thus the point 
\begin{equation}\label{eq:82}
C_n = A_n - \tau A_0
\end{equation}
is not only an invariant point but also this point is 
intrinsically connected with the geometry of a lightlike 
hypersurface $V^{n-1}$ and is defined by elements of a 
fourth-order neighborhood of a point $x \in V^{n-1}$. 

We proved the following result.

\begin{theorem}
On a lightlike hypersurface $V^{n-1}$ of a manifold $(M, c)$,  
 a complete invariant normalization intrinsically connected 
 with the geometry of  $V^{n-1}$ is defined by elements of a 
fourth-order differential neighborhood of a point 
$x \in V^{n-1}$. This complete normalization induces an affine 
connections that is also intrinsically connected with the 
geometry of  $V^{n-1}$.  
\end{theorem} 

The last statement of Theorem 10 follows from Theorem 8. 
The normalizing hyperplane $L_x$ discussed in Theorem 8 
is defined by the harmonic pole $H$ of the point $x$ 
with respect to the singular points of the isotropic 
geodesic $A_0 A_1$ of the hypersurface $V^{n-1}$; the normalizing 
$(n-3)$-plane $R_x$  which is the intersection of 
 the screen subspace and the complementary screen subspace 
(both belong to the hyperplane $T_x (V^{n-1}$)); and 
finally the invariant point $C_n$ (defined by formula (82) on 
the normalizing straight line $A_0 A_n$) which is 
conjugate to the screen subspace $S_x$ with respect to the 
isotropic cone $C_x$. 

{\bf 2.} Let us clarify a geometric meaning 
of the normalizing point $C_n$ on the 
straight line $A_0 A_n$. We consider 
this straight line as the line belonging to 
a local hyperquadric $(Q^n_1)_x$ that is tangent  to a 
manifold $(M, c)$ at a point $x$.  
On the manifold $(M, c)$, to this line $A_0 A_n$ 
 there corresponds an   isotropic geodesic $\widetilde{l}$ 
which we will also denote by $A_0 A_n, \;\; \widetilde{l} 
= A_0 A_n$.

When a point $A_0$ describes a lightlike hypersurface 
$V^{n-1} \subset (M, c)$, the   isotropic geodesic $\widetilde{l}$ 
describes an {\em isotropic congruence} $U$. As a point manifold, 
this congruence is an $n$-dimensional domain 
on the manifold $(M, c)$. This is the reason that 
we will denote it by $U^n$. Moreover, $U^n = \widetilde{f} 
(V^{n-1} \times \widetilde{l})$, where 
$\widetilde{f}$ is a differentiable mapping of the direct 
product $V^{n-1} \times \widetilde{l}$ onto $(M, c)$.

Let us find the Jacobian of the mapping $\widetilde{f}$. 
To this end we consider an arbitrary point $Z = A_n + z A_0$ 
on the    isotropic geodesic $\widetilde{l}$. The 
differential of this point has the form
$$
\renewcommand{\arraystretch}{1.3}
\begin{array}{ll}
d Z =\!\!\!\!& - \omega^1_1 Z + [dz + z (\omega_0^0 + \omega_1^1) 
+ \omega_n^0] A_0 \\
\!\!\!\!& + (\tau_b^a + z \delta_b^a) \omega^b A_a 
+ (z A_1 + \tau^a A_a - A_{n+1}) \omega^1,
\end{array}
 \renewcommand{\arraystretch}{1.3}
$$
where $\tau^a = g^{ab} \tau_b$ and $\tau_a^b = g^{ac} \tau_{cb}$. 
At regular points $Z$ of the mapping $\widetilde{f}$. the tangent 
space to the domain $U^n$ coincides with the $n$-dimensional 
tangent subspace of the manifold $(M, c)$ which is defined by 
the points $Z, A_0, A_n$, and $A_{n+1} - z A_1 - \tau^a A_a$. 
At singular points $Z$, the dimension of this tangent space 
is reduced. This happened at the points at which the rank of 
the Jacobian matrix of the mapping $\widetilde{f}$,
$$
\widetilde{J} = (\tau_b^a + z \delta_b^a), \;\;\; a, b = 
2, \ldots, n-1,
$$ 
is reduced, and only at such points. At such points 
the Jacobian vanishes:
\begin{equation}\label{eq:83}
\det  (\tau_b^a + z \delta_b^a) = 0.
\end{equation}

Denote the eigenvalues of the matrix $(\tau_b^a)$ by 
$\widetilde{\tau}_a$. In general, the matrix $(\tau_b^a)$ 
is not symmetric, since the matrix $(\tau_{ab})$ occurring 
in equation (77) is not symmetric. Thus 
the eigenvalues $\widetilde{\tau}_a$ of the matrix $(\tau_b^a)$ 
can be either real or complex conjugate. Hence the solutions 
of equation (81), that are defined by the 
formula $z_a = - \widetilde{\tau}_a$, as well 
as the singular points 
\begin{equation}\label{eq:84}
Z_a = A_n -\widetilde{\tau}_a A_0
\end{equation}
 can be also complex conjugate.

But even in the case of complex conjugate roots of equation (83), 
by Vieta's theorem, the sum of the roots of equation (83) is 
the negative trace of the matrix $(\tau_b^a)$, \begin{equation}\label{eq:85}
\sum_{a=2}^{n-1} z_a = - \sum_{a=2}^{n-1} \tau_a = 
= - \sum_{a=2}^{n-1} \tau_b^b = - \tau_{ab} g^{ab}.
\end{equation}
This relation allows us to find a geometric meaning of 
the invariant point $C_n$ on the isotropic geodesic 
$\widetilde{l} = A_0 A_n$. Since the coefficient $\tau$ 
occurring in the formula (82) defining the point $C_n$ is determined 
by equation (79), {\em the point $C_n$ is the 
harmonic pole of the point $A_0$ with respect to the singular 
points $Z_a$ of the mapping $\widetilde{f}$.}

\section{Integrability of screen distributions}

 Consider the screen distribution $S$ 
and  the  complementary  screen distribution $\widetilde{S}$ 
of a lightlike hypersurface $V^{n-1}$ 
which we have constructed in Section {\bf 5.4}. 
The distribution $S$ is formed by the subspaces 
$S_x = A_0 \wedge  C_2 \wedge \ldots \wedge C_{n-1}$, where the 
points $C_a$ are defined by equations (75), and the 
distribution $\widetilde{S}$ is formed by the subspaces 
$\widetilde{S}_x = H \wedge  C_2 \wedge \ldots \wedge C_{n-1}$, 
where $H$ is the harmonic pole of the point $A_0$ with respect to 
the singular points of the isotropic geodesic $l = A_0 A_1$ 
of the hypersurface $V^{n-1}$. 

As we indicated in Section {\bf 6.1}, assuming that the tensors 
$h_b^a$ and $H_b^a$ are nondegenerate, we can reduce a frame 
bundle associated with the hypersurface $V^{n-1}$ in such a way 
that the points $A_1$ and $H$ as well as the points $A_a$ and $C_a$ 
will be superposed. This gives 
$$
S_x = A_0 \wedge  A_2 \wedge \ldots \wedge A_{n-1}, \;\;
\widetilde{S}_x = A_1 \wedge  A_2 \wedge \ldots \wedge A_{n-1}.
$$
As a result, the screen distribution $S$ of 
the hypersurface $V^{n-1}$ is defined by the 
differential equation
\begin{equation}\label{eq:86}
\omega_0^1 = 0,
\end{equation}
and   the  complementary  screen distribution $\widetilde{S}$ is 
defined by the equation 
\begin{equation}\label{eq:87}
\omega_1^0 = 0.
\end{equation}

Let us find the conditions of integrability of the screen 
distributions $S$ and $\widetilde{S}$. Since in the reduced 
frame bundle the geometric object $\mu_a$ vanishes, 
equation (55) connecting the 1-forms $\omega_1^0$ 
and $\omega_0^1$ takes the form
\begin{equation}\label{eq:88}
\omega_1^0 = \mu \omega_0^1,
\end{equation}
where the invariant $\mu$ is different from 0. 
Thus {\em if one of the distributions $S$ and $\widetilde{S}$ 
defined by equations $(86)$ and $(87)$, respectively, 
is integrable, then another one is also integrable}. 

From the structure equations (2)---(5) of a manifold $(M, c)$ 
endowed with the pseudoconformal structure $CO (n-1, 1)$ 
it follows that 
$$
\renewcommand{\arraystretch}{1.3}
\begin{array}{ll}
d\omega_0^1 = \omega_0^0 \wedge \omega_0^1 + \omega_0^1 \wedge 
\omega_1^1 + \omega_0^a \wedge \omega_a^1, \\
d\omega_1^0 = \omega_1^0 \wedge \omega_0^0 + \omega_1^1 \wedge 
\omega_1^0 + \omega_1^a \wedge \omega_a^0 + 2 C_{11a} 
\omega^1 \wedge \omega^a + C_{1ab} \omega^a \wedge \omega^b.
\end{array}
\renewcommand{\arraystretch}{1}
$$
Substituting 
the values of 
the 1-forms $\omega_a^1$ and 
$\omega_a^0$ taken from equations (77) into the last two equations, 
we find that 
$$
\renewcommand{\arraystretch}{1.3}
\begin{array}{ll}
d\omega_0^1 \equiv \tau_{ab}  \omega^a \wedge \omega^b \pmod{\omega_0^1}, \\
d\omega_1^0 \equiv (h_a^c \sigma_{cb} + C_{1ab}) 
  \omega^a \wedge \omega^b \pmod{\omega_1^0}.
\end{array}
\renewcommand{\arraystretch}{1}
$$
This and the Frobenius theorem (see, for example, 
[BCGGG 91]) imply that the condition of integrability 
of the screen distribution $S$ has the form 
\begin{equation}\label{eq:89}
\tau_{ab} = \tau_{ba}, 
\end{equation}
and the condition of integrability 
of the complementary screen distribution $\widetilde{S}$ 
has the form 
\begin{equation}\label{eq:90}
h_a^c \sigma_{cb} + C_{1ab} = h_b^c \sigma_{ca} + C_{1ba}.
\end{equation}

Let us prove that conditions (89) and (90) are equivalent. 
In fact, in the reduced frame bundle we have equations (76) 
and (77). Substituting 
the values of $\omega_a^0$ and 
$\omega_a^1$ taken from (77) into the first equation of (76) 
and using the fact that the forms $\omega^1$ and $\omega^a, \;
a = 2, \ldots , n - 1$, are linearly independent, 
we find that 
$$
- h_a^b \sigma_b - \mu \tau_a + 2 C_{11a} = 0; \;\;
- h_a^c \sigma_{cb} - \mu \tau_{ab} = C_{1ab} + \nu_{ab}.
$$
Alternating the second equation with respect to the indices 
$a$ and $b$, we obtain that 
\begin{equation}\label{eq:91}
h_a^c \sigma_{cb} - h_b^c \sigma_{ca}  + C_{1ab} - C_{1ba} = 
- 2\mu \tau_{[ab]}.
\end{equation}
Since $\mu \neq 0$, it follows from (91) that conditions 
(89) and (90) are equivalent.

Suppose that the screen distribution $S$ of a lightlike 
hypersurface $V^{n-1}$ is integrable. Denote by $S (x)$ an 
integral manifold of this distribution passing through 
a point $x = A_0$ of the isotropic geodesic $l$ of 
the hypersurface $V^{n-1}$. Then Pfaffian 
equation (86) determines 
also a stratification of the normalizing congruence $U^n$ of 
the hypersurface $V^{n-1}$ into a one-parameter family 
of lightlike hypersurfaces $U^{n-1} (x)$ whose generators 
$\widetilde{l} = A_0 A_n$ pass through the points of the manifold 
$S (x)$.

The following theorem combines the results obtained in 
this subjection.

\begin{theorem}
If on a lightlike hypersurface $V^{n-1}$ the conditions 
$(89)$ hold, then
\begin{description}
\item[(a)] The screen distribution 
$S$ and the complementary screen distribution $\widetilde{S}$, 
are integrable.

\item[(b)] The normalizing lightlike congruence $U^n$ 
of a hypersurface $V^{n-1}$ is stratified into a one-parameter 
family of lightlike hypersurfaces $U^{n-1}$.

\item[(c)] All singular points of the congruence $U^n$ 
are real and coincide with singular points 
of the hypersurfaces $U^{n-1}$. 
\end{description}
\end{theorem}

Part (c) of Theorem 11 follows from the symmetry 
of the matrix $(\tau_{ab})$: all eigenvalues of such a matrix 
(roots of (83)) are real.

Note also that if a $CO (n-1, 1)$-structure on 
the manifold 
$(M, c)$ is conformally flat, then equations (76) and (89) imply 
that the affinors $\tau_b^a$ and $\sigma_b^a = g^{ac} \sigma_{cb}$ 
of a lightlike hypersurface $V^{n-1}$ 
are diagonalized simultaneously with the affinors $h_b^a$ and $H_b^a$. 
Geometrically this means that the torses 
formed by isotropic geodesics on every lightlike hypersurface $U^{n-1}$ 
correspond one to another. 

\section{Construction of affine connections intrinsically 
connected with  a lightlike hypersurface}

{\bf 1}. In Section {\bf 4}, we have already considered 
 the question of finding 
an affine connection on a lightlike hypersurface $V^{n-1}$ 
of a manifold $(M, c)$.  As we proved in Theorem 10, 
an  affine connection on $V^{n-1}$---denote it by $\gamma_1$---is defined 
in a fourth-order differential neighborhood. Equations (47) show 
that this connection is torsion-free. For finding the curvature 
tensor of the connection $\gamma_1$ we substitute the values (77) 
of the forms 
$\omega_a^0$ and $\omega_a^1$ into (48). 
In addition, since the vertex $A_n$ coincides with the 
harmonic pole $C_n$ of the point $A_0 = x$ with respect 
to the singular points $Z_a$ of the isotropic geodesic $A_0 A_n$, 
the 1-form $\omega_n^0$ 
becomes a principal form:
\begin{equation}\label{eq:92}
\omega_n^0 = \varphi_1 \omega^1 + \varphi_a \omega^a.
\end{equation}
The coefficients $\varphi_1$ and $\varphi_a$ in (92) are defined 
in a fifth-order neighborhood of a point $x \in V^{n-1}$. 
This implies that the curvature tensor of 
the  affine connection $\gamma_1$ induced by the invariant 
normalization of $V^{n-1}$ we have constructed is also defined 
in a fifth-order neighborhood of a point $x \in V^{n-1}$.

However, there is another way to construct an affine connection 
on a hypersurface $V^{n-1}$. To this end, we note that in a 
reduced frame bundle associated with 
 a third-order neighborhood of a point $x \in V^{n-1}$ 
the 1-forms $\omega_1^a$ and $\omega_a^1$ in equations (47) 
become principal: they are expressed by formulas (42) and (77). 
Substituting the expressions of these forms into equations 
(47), we find that 
\begin{equation}\label{eq:93}
\renewcommand{\arraystretch}{1.3}
\left\{
\begin{array}{ll}
d \omega^1 = \omega^1 \wedge (\omega_1^1 - \omega^0_0) 
+ \omega^a \wedge (\tau_a \omega^1 + \tau_{ab} \omega^b), \\
d \omega^a = \omega^b 
\wedge (\omega_b^a  - \delta_b^a \omega^0_0) 
+ h_b^a \omega^1 \wedge \omega^b.
\end{array} 
\right.
\renewcommand{\arraystretch}{1}
\end{equation}
Thus the forms $\omega_1^1 - \omega^0_0$ and $\omega_b^a  - 
\delta_b^a \omega^0_0$ can be considered as the only  connection 
forms of the affine connection $\gamma_2$, and the 2-forms 
\begin{equation}\label{eq:94}
\renewcommand{\arraystretch}{1.3}
\left
\{\begin{array}{ll}
\Theta^1 
  =  \omega^a \wedge (\tau_a \omega^1 + \tau_{ab} \omega^b), \\
\Theta^a = h_b^a \omega^1 \wedge \omega^b.
\end{array} 
\right.
\renewcommand{\arraystretch}{1}
\end{equation}
as the torsion 2-forms of this connection. 

In the decompositions (48) of exterior differentials of 
the 1-forms $\omega_1^a$ and $\omega_a^1$ we have only the forms 
$\omega_1^0, \; 
\omega_a^0$, and $\omega_n^a = g^{ab} \omega_1^a$, 
and the form $\omega_n^0$ does not occur in these 
decompositions. Thus {\em the torsion and 
curvature tensors of the affine 
connection $\gamma_2$ are defined in 
 a fourth-order neighborhood of a point $x \in V^{n-1}$}, not 
the fifth-order as  was the case for the connection $\gamma_1$. 

If the principal parameters are fixed, then the first and the last 
subsystems of system (48) take the form
\begin{equation}\label{eq:95}
\renewcommand{\arraystretch}{1.3}
\left\{
\begin{array}{ll}
\delta (\pi_1^1 - \pi^0_0) = 0, \\
\delta  (\pi_b^a - \delta_b^a \pi^0_0)  
=  (\pi_b^c - \delta_b^c \pi_0^0) 
\wedge (\pi_c^a  - \delta_c^a \pi_0^0).
\end{array} 
\right.
\renewcommand{\arraystretch}{1}
\end{equation}
Equations (95) are the structure equations of the group $G_2$ 
defining the connection $\gamma_2$ on 
the hypersurface $V^{n-1}$. It follows from (95) that 
the group $G_2$ is the direct product of the group 
${\bf R}^+$ of homotheties (${\bf R}^+$ is the 
multiplicative group of positive real  numbers), 
and the general 
linear group $GL (n-2, {\bf R})$ over the field of real numbers: 
$$
G_2 = {\bf R}^+ \times {\bf GL} (n-2, {\bf R}).
$$

It follows from equations (48) that the curvature forms 
of the connection $\gamma_2$ can be written as 
$$
\left\{
\begin{array}{ll}
\Omega^1_1 =  \omega^a_0 \wedge \omega^0_a + \omega_1^a  
\wedge \omega_a^1 + C^1_{1kl} \omega^k \wedge \omega^l, \\
\Omega^a_b =  \omega_b^0 \wedge \omega_0^a + \omega_b^1 \wedge 
\omega_1^a + \omega_b^n \wedge \omega_n^a + g_{bc} g^{ae} 
\omega_0^c \wedge \omega_e^0 
- \delta_b^a \omega_0^c \wedge \omega_c^0 + C_{bkl}^a \omega^k 
\wedge \omega^l,
\end{array} 
\right.
\renewcommand{\arraystretch}{1}
$$
where $a, b = 2, \ldots , n-1; \; k, l = 1, 2, \ldots , n-1$. 
Substituting the values (77) of the forms 
$\omega_a^0, \; \omega_a^1$, and $\omega_n^a = g^{ab} \omega_b^1$ 
into these expressions and using (42), we find the following 
values of the curvature 2-forms $\Omega_1^1$ and $\Omega_b^a$:
\begin{equation}\label{eq:96}
\left\{
\begin{array}{ll}
\Omega^1_1 = & \!\!\!\! (2 C^1_{11a} - \sigma_a - h_a^c \tau_c) 
\omega^1 \wedge \omega^a + (\sigma_{ab} + h_a^c \tau_{cb} + C^1_{1ab}) 
\omega^a  \wedge \omega^b,  \\
\Omega^a_b = & \!\!\!\!  
(\delta_b^a \sigma_c + \delta^a_c \sigma_b 
- g_{bc} \sigma^a + h_c^a \tau_b - h_{bc} \tau^a 
+ 2 C^a_{b1c}) \omega^1 \wedge \omega^c, \\
&  \!\!\!\!  
+ (\delta_e^a \sigma_{bc} + g_{bc} \sigma^a_e + h_e^a \tau_{bc} 
 + h_{bc} \tau^a_e - \delta_b^a \sigma_{ce} 
+  C^a_{bce}) \omega^c  \wedge \omega^e,
\end{array} 
\right.
\renewcommand{\arraystretch}{1}
\end{equation}
where $\sigma^a = g^{ab} \sigma_b, \; \tau^a = g^{ab} \tau_b, \; 
\sigma_b^a = g^{ac} \sigma_{cb}$, 
and $\tau_b^a = g^{ac} \tau_{cb}$. It follows that the components 
of the curvature tensor of the connection $\gamma_2$ 
are determined by the following formulas:
\begin{equation}\label{eq:97}
\left\{
\begin{array}{ll}
R^1_{11a} = & \!\!\!\! 2 C^1_{11a} - \sigma_a - h_a^c \tau_c, \\
R^1_{1ab} = &\!\!\!\! 
\sigma_{[ab]} + h_{[a}^c \tau_{|c|b]} + C^1_{1ab},\\ 
R^a_{b1c} = &\!\!\!\! \delta_b^a \sigma_c + \delta^a_c \sigma_b 
- g_{bc} \sigma^a + h_c^a \tau_b - h_{bc} \tau^a 
+ 2 C^a_{b1c}, \\
R^a_{bce} = &\!\!\!\!  \delta_{[e}^a \sigma_{|b|c]} + g_{[b|c|} 
\sigma^a_{e]} + h_{[e}^a \tau_{|b|c]} 
 + h_{b[c} \tau^a_{e]} - \delta_b^a \sigma_{[ce]} 
+  C^a_{bce}.
\end{array} 
\right.
\renewcommand{\arraystretch}{1}
\end{equation}

\noindent
{\em Authors' addresses}:\\

\noindent
\begin{tabular}{ll}
M.A. Akivis &                           V.V. Goldberg \\
Department of Mathematics      &      Department of Mathematics\\
Jerusalem College of Technology  & 
                       New Jersey Institute of Technology \\
-- Mahon Lev, P. O. B. 16031 & University Heights \\
Jerusalem 91160, Israel & Newark, NJ 07102, U.S.A. \\
\\
E-mail address: akivis@avoda.jct.ac.il & E-mail address: 
                                        vlgold@numerics.njit.edu    
\end{tabular}

\end{document}